\let\orilabel\label
\documentclass[review,hidelinks,onefignum,onetabnum]{siamart250211}
\let\label\orilabel

\usepackage{lipsum}
\usepackage{amsfonts}
\usepackage{graphicx}
\usepackage{epstopdf}
\usepackage{algorithmic}
\usepackage{multirow}
\usepackage{tikz}
\usepackage{makecell}
\usetikzlibrary{patterns}
\ifpdf
  \DeclareGraphicsExtensions{.eps,.pdf,.png,.jpg}
\else
  \DeclareGraphicsExtensions{.eps}
\fi
\usepackage{comment}

\newsiamremark{remark}{Remark}
\newsiamremark{hypothesis}{Hypothesis}
\crefname{hypothesis}{Hypothesis}{Hypotheses}
\newsiamthm{claim}{Claim}
\newsiamremark{fact}{Fact}
\crefname{fact}{Fact}{Facts}

\newcommand{\nfa}{\textrm{NFA}}

\newcommand{\bn}{\mathbf{n}}
\newcommand{\bv}{\boldsymbol{v}}

\newcommand{\sett}[2]{\left\{{#1}\,:~{#2}\right\}}

\newcommand{\sym}{\nabla_{sym}}

\newcommand{\gd}{\nabla}

\DeclareMathOperator{\dv}{div}
\usepackage{xcolor}

\definecolor{coolwarmblue}{HTML}{3B4CC0}
\definecolor{coolwarmred}{HTML}{B40426}

\headers{Modeling and Simulation of Performance in OPV}{P. \c{C}\.{i}lo\u{g}lu et. al}

\title{Modeling and Simulation of Device Performance in Organic Photovoltaics\thanks{Submitted to the editors.}}

\author{
Pel\.{i}n \c{C}\.{i}lo\u{g}lu$^{\diamond,}$\thanks{Faculty of Mathematics, University of Technology Chemnitz, 09107, Chemnitz, Germany   (\email{pelin.ciloglu@mathematik.tu-chemnitz.de}, \email{martin.stoll@mathematik.tu-chemnitz.de}).}
\and Carmen Tretmans$^{\diamond,}$\thanks{Faculty of Mathematics, University of Augsburg, 86159, Augsburg, Germany  (\email{carmen.tretmans@uni-a.de},  \email{jan-f.pietschmann@uni-a.de}).}
\and Carsten Deibel\thanks{Institute of Physics, University of Technology Chemnitz, 09107, Chemnitz, Germany (\email{deibel@physik.tu-chemnitz.de}).}
\and Roderick MacKenzie\thanks{Department of Engineering, Durham University, Durham, United Kingdom (\email{roderick.mackenzie@durham.ac.uk}).}
\and
Roland Herzog\thanks{Interdisciplinary Center for Scientific Computing (IWR), Heidelberg University, 69120, Heidelberg, Germany (\email{roland.herzog@iwr.uni-heidelberg.de}).} 
\and Jan-F. Pietschmann\footnotemark[3]\;\,$^{,}$\thanks{Centre for Advanced Analytics and Predictive Sciences (CAAPS), University of Augsburg,
Universit\"{a}tsstr. 12a, 86159 Augsburg, Germany}
\and
Martin Stoll\footnotemark[2]
}

\usepackage{amsopn}

\begin{document}
\nolinenumbers
\maketitle
\def\thefootnote{$\diamond$}\footnotetext{These authors contributed equally to this work.}

\begin{abstract}
We present a pipeline to study the device performance of organic solar cells in silico. We introduce a mathematical model that includes the dynamics of excitons as well as their dissociation at bulk heterojunctions within the nanomorphology of the active layer. This is combined with realistic morphologies that we obtain from a detailed phase field model. To solve the coupled nonlinear system, we use a finite element discretization, robust linear solvers, and three numerical schemes, Newton, Gummel, and Semi--Newton--Gummel. This allows for an efficient simulation of the complete OPV device and results in current-voltage curves that can readily be compared to measured data.
\end{abstract}



\section{Introduction}
Organic photovoltaics (OPV) have been shown to be a promising technology, with potential for being a low-cost, sustainable, and flexible alternative to inorganic photovoltaics \cite{Wang2023, Wopke2022}. While such inorganic cells (e.g., silicon-based) still outperform their organic counterparts in terms of efficiency and stability, the development of OPVs has made great progress in recent years. For instance, the introduction of non-fullerene acceptors (NFA), having highly tunable photoelectric properties, within the active layer of the cell considerably improved efficiency and stability \cite{Zhang_2022}. This paper introduces a mathematical model for charge generation and transport, taking the specifics of organic based materials into account. Subsequently, a numerical scheme allowing for its efficient simulations is developed. This results in a pipeline allowing to examine the electrical performance in silico.

The production of OPVs is based on a solvent-based fabrication process. A blend of polymer acting as the electron donor, non-fullerene acceptor (NFA), and solvent is deposited on a substrate. While polymer and NFA undergo spontaneous phase separation, the solvent will evaporate, drying the thin film and eventually halting the phase separation process. The resulting thin film consists of regions of pure donor and acceptor materials with an irregular self-assembled nanostructure, as shown schematically in Figure~\ref{fig_sketchcell}. This generated morphology of the active layer is classified by so-called bulk heterojunctions (BHJ), indicating the irregular donor--acceptor interfaces within the bulk. This is in contrast to most inorganic cells, containing bilayer junctions where donor and acceptor are stacked. 

The efficiency of OPVs depends not only on the generated morphology in the active layer, but also on material properties of donor and acceptor, like the energy levels of the HOMO (\textit{Highest Occupied Molecular Orbital}) and LUMO (\textit{Lowest Unoccupied Molecular Orbital}), see \cite{burke2015device} for details. Although organic and inorganic photovoltaics share the general idea of converting incident photons to electrical power, the underlying mechanisms governing charge generation and transport differ significantly, requiring distinct modeling approaches. In traditional, inorganic photovoltaics, incident photons directly generate free charge carriers within its crystal lattice. Charges are then transported via continuous band transport. In contrast, free charge carriers cannot be directly generated within organic cells. Incident photons create bound electron-hole pairs, so-called \textit{excitons}. The Coulomb attraction between the hole and the electron is sufficiently strong such that the thermal energy is insufficient to dissociate the exciton immediately. To split the exciton into free charges, an additional driving force is needed, which is found in the energy offset of donor and acceptor materials at the bulk heterojunction, indicating the impact of material energy levels on device performance. Consequently, for the exciton to be dissociated, it is required to reach an acceptor-donor interface. Naturally, the exact active layer morphology has a huge impact on the performance of the organic cell, too. This interplay of factors shows the necessity of a tailored electrical model describing OPVs that include material properties as well as active layer morphologies. 

Various models have been proposed to study the electrical performance of OPVs, which extend classical drift--diffusion models used in the inorganic case \cite{markowich1990semiconductor}. In \cite{Brinkman2013modell_asymptotics}, a drift--diffusion--reaction model for organic bilayers is presented, including the process of exciton generation and dissociation. Models for the organic--specific hopping behavior of charges within the active material are proposed in \cite{Doan2019, DFF20, Glitzky2024stack}, presenting models for charge carrier transport in organic semiconductor devices.

In recent decades, various numerical methods have been introduced for solving semiconductor device equations. These equations are strongly coupled and highly nonlinear, resulting in significant challenges for the numerical schemes \cite{markowich1990semiconductor, Markowich_1985}. The most common discretization techniques used in semiconductor simulations are the finite volume method (FVM), the finite element method (FEM), and the discontinuous Galerkin method (DG). The finite volume scheme combined with the Scharfetter–Gummel discretization \cite{Scharfetter_Gummel_1969} of the drift–diffusion equations, commonly referred to as FVSG \cite{Bank_Rose_Fichtner_1983}, has emerged as the standard numerical strategy in industrial semiconductor simulators. Although the finite volume method ensures local conservation and is robust for drift–diffusion systems, it may be less accurate on unstructured meshes and difficult to extend to higher-order approximations. DG schemes  \cite{cockburn2012discontinuous, cockburn2006introduction} such as hybridizable discontinuous Galerkin (HDG) schemes \cite{cockburn2009unified} are particularly well-known for semiconductor simulations due to their local conservation, stability, and high-order accuracy \cite{Brinkman2013modell_asymptotics, chen2019hdg}. 
The classical finite element method \cite{brenner2008mathematical, strang1973analysis} is fundamental to the numerical solution of the drift–diffusion model \cite{Bank_Rose_Fichtner_1983, barnes1977finite, buturla1981finite}, as it supports higher-order approximations and can be applied to domains of arbitrary geometric complexity. In contrast to DG methods, the number of degrees of freedom remains at a manageable level, and unlike finite volume schemes, FEM naturally handles complex geometries. For these reasons, the finite element approach is adopted in this work. 

The resulting nonlinear system is treated using three different strategies: the classical Newton method, the standard Gummel iteration, and a hybrid Semi--Newton--Gummel approach. The Newton method is known for its fast convergence, but constructing and solving the full Jacobian system can be computationally demanding. In contrast, the Gummel iteration \cite{Gummel_1964} decouples the equations and solves them sequentially and semi-implicitly, significantly reducing the computational effort. The hybrid strategy balances the strengths of both approaches by applying Newton-type updates to the individual subsystems while maintaining the overall decoupled structure of the Gummel approach. This combination provides a robust and efficient alternative, especially in situations where the strong nonlinear coupling characteristic of drift–diffusion models presents challenges for fully coupled Newton solvers.

The aim of this paper is to predict the electrical performance of OPVs based on an extended model that uses realistic morphologies of the active layer, all while incorporating the organic device specific exciton dynamics. 
We make the following contributions to existing models:
\begin{itemize}
    \item Explicit modeling of exciton dynamics within nanomorphologies.
    \item Coupling of charge generation by exciton splitting to realistic active layer morphologies.
    \item Inclusion of a morphology-dependent charge transport to the electrodes.
    \item Development of an efficient hybrid numerical approach allowing for simulations of realistic morphology-dependent current-density--voltage characteristics ($I$-$V$ characteristics). 
    \item Case studies in up to three spatial dimensions.
\end{itemize}

\section{Mathematical model}
In this section, we will derive the governing equations for electron and hole densities $n$ and $p$, the exciton density $X$, and the electric potential $\psi$. To incorporate the realistic nanostructure of the active layer, we employ the morphology evolution governed by a phase field model described in our previous work \cite{PCiloglu_2025}. We first introduce the domain of our model. Then, we define the equation that governs the behavior of the excitons, necessary for the generation of free charge carriers. With this in place, we continue with the charge transport model and the electrical potential. 
\subsection{Domain}
Since our main focus is the influence of the active layer morphology on the electric device performance, we restrict the domain to cover solely the organic active layer $\Omega_{\mathrm{org}} \in \mathbb{R}^d$, $d=2,3$, thereby neglecting the charge behavior within the two electrodes. To describe the morphology of the active layer, we make use of a continuous phase field $\phi: \Omega_{\mathrm{org}} \to [0,1]$ describing the volume fraction of acceptor material. Thus, $\phi = 1$ represents areas of pure acceptor material, $\phi = 0$ pure donor, and interfaces correspond to values $0<\phi<1$. The interfaces between donor and acceptor regions within the organic domain can thus be located by evaluating the phasefield gradient $\lambda |\nabla \phi|$, where $\lambda$ represents the characteristic interface thickness. While the process of morphology formation is certainly of interest, here, we solely focus on the electrical performance for given morphologies. However, since we still wish to keep the simulations as realistic as possible, we employ simulation results on realistic device morphologies presented in our previous work \cite{PCiloglu_2025}. 

To introduce the various boundary conditions, we denote the top and bottom boundary, representing the interface with the top and bottom electrode by $\Gamma_{\mathrm{top/org}}$ and $\Gamma_{\mathrm{bot/org}}$, respectively. The remaining boundary part represents the insulated part of the active layer and is denoted by $\Gamma_{\mathrm{ins}}$. A sketch of the geometry can be found in Figure~\ref{fig_sketchdomain}. In addition, we fix a maximal simulation time $t_{\mathrm{max}} > 0$.
As already discussed, the energy offset between donor (polymer) and acceptor (NFA) material is crucial for charge separation. We define the constant, material dependent HOMO and LUMO levels by $E^i_H$ and $E^i_L$, $i \in \{p, \nfa \}$, respectively, and employ a linear interpolation for the pointwise HOMO and LUMO levels over the organic domain, $E_{H/L} = E_{H/L}^{p} + (E_{H/L}^{\nfa} - E_{H/L}^{p}) \phi$. We furthermore define the constant temperature $T$, the elementary charge $q_e$, and the Boltzmann constant $k_B$. 

\begin{figure}[htp!]
  \centering
  \begin{minipage}{0.44\textwidth}
  \centering
    \begin{tikzpicture}[scale=4]

  \def\hel{0.15} 
  \def\hactive{0.7} 

  \fill[pattern=north west lines, pattern color=gray!50] (0,0) rectangle (1,\hel);
  \node at (0.5, \hel/2) {bottom electrode};

  \begin{scope}
      \clip (0,\hel) rectangle (1, 1-\hel);

      
      \begin{scope}[shift={(-0.433, \hel)}, scale={\hactive/3}]
          
         \fill[coolwarmblue!70] (0,0) rectangle (8,3);

         \fill[coolwarmred!70] (0.7, 3) .. controls (0.7, 2.7) and (1.0, 2.3) .. (1.3, 2.0)
         .. controls (1.5, 1.8) and (1.1, 1.5) .. (1.4, 1.2)
         .. controls (1.6, 1.0) and (1.2, 0.8) .. (1.4, 0.5)
         .. controls (1.6, 0.3) and (1.8, 0.5) .. (2.0, 0.7)
         .. controls (2.2, 1.0) and (1.9, 1.5) .. (2.1, 1.8)
         .. controls (2.3, 2.1) and (1.8, 2.5) .. (2.1, 3)
         -- cycle;

         \begin{scope}[yscale=-1,xscale=1, shift = ({0,-3})]
             \fill[coolwarmred!70] (2.4, 3)
             .. controls (2.3, 2.5) and (2.5, 2.1) .. (2.7, 1.8)
             .. controls (2.9, 1.5) and (2.6, 1.3) .. (2.9, 1.0)
             .. controls (3.1, 0.6) and (2.7, 0.4) .. (2.8, 0.0)
             -- (3.1, 0.0)
             .. controls (3.3, 0.4) and (3.4, 0.8) .. (3.6, 1.1)
             .. controls (3.7, 1.4) and (3.4, 1.7) .. (3.6, 2.1)
             .. controls (3.3, 2.4) and (3.0, 2.8) .. (3.0, 3.0)
             -- cycle;
         \end{scope}

         \fill[coolwarmred!70] (4.3, 3)
         .. controls (4.3, 2.8) and (4.5, 2.5) .. (4.7, 2.1)
         .. controls (4.9, 1.8) and (4.5, 1.5) .. (4.8, 1.2)
         .. controls (5.0, 1.0) and (4.6, 0.8) .. (4.9, 0.5)
         .. controls (5.1, 0.3) and (4.8, 0.2) .. (4.6, 0.3)
         .. controls (4.4, 0.5) and (4.3, 1.0) .. (4.2, 1.6)
         .. controls (4.1, 2.0) and (3.9, 2.5) .. (4.1, 3)
         -- cycle;
         
         \begin{scope}[yscale=-1,xscale=1, shift = ({0,-3})]
             \fill[coolwarmred!70] (5.0, 3)
             .. controls (5.05, 2.7) and (5.25, 2.3) .. (5.55, 2.0)
             .. controls (5.35, 1.6) and (5.65, 1.3) .. (5.45, 1.0)
             .. controls (5.55, 0.7) and (5.25, 0.5) .. (5.45, 0.3)
             .. controls (5.65, 0.1) and (5.85, 0.3) .. (6.05, 0.5)
             .. controls (6.15, 0.7) and (5.85, 1.4) .. (6.05, 1.5)
             .. controls (6.15, 1.7) and (5.85, 2.1) .. (6.0, 2.5)
             .. controls (6.05, 2.6) and (5.85, 2.9) .. (5.75, 3)
             -- cycle;
       \end{scope}

        \fill[coolwarmred!70] (6.3, 3)
         .. controls (6.2, 2.5) and (6.5, 2.1) .. (6.7, 1.6)
         .. controls (6.1, 1.2) and (7.1, 1.5) .. (6.8, 0.5)
         .. controls (6.5, 0.4) and (6.8, 0.25) .. (6.7, 0.)
         -- (7.1, 0.0)
         .. controls (7.4, 0.2) and (7.3, 0.9) .. (7.4, 1.3)
         .. controls (7.3, 1.5) and (6.9, 2.0) .. (7.1, 2.2)
         .. controls (7.4, 2.1) and (7.5, 1.7) .. (7.7, 1.5)
         .. controls (7.8, 1.7) and (7.4, 2.3) .. (7.5, 2.7)
         .. controls (7.5, 2.8) and (7.2, 3.0) .. (7.0, 3)
         -- cycle;

      \end{scope}
  \end{scope}

  \node[fill=white, fill opacity=0.7, text opacity=1, rounded corners=2pt, inner sep=2pt] 
       at (0.5, 0.5) {organic active layer};

  \fill[pattern=north west lines, pattern color=gray!50] (0,1-\hel) rectangle (1,1);
  \node at (0.5, 1-\hel/2) {top electrode};

  \draw (0,0) rectangle (1,1);
\end{tikzpicture}
\label{fig_sketchcell}
\caption{Sketch of an organic solar cell containing an organic active layer, top, and bottom electrode. The active layer consists of self-assembled pure donor (polymer) and acceptor (NFA) regions, indicated by the blue and red areas.}
  \end{minipage} \hspace{0.75cm} 
  \begin{minipage}{0.44\textwidth} 
    \centering
    \vspace{-0.2cm}
    \begin{tikzpicture}[scale=4]

  \def\hel{0.15} 

    \draw[dashed] (0,0) rectangle (1,1);
    \node at (0.5, \hel/2) {bottom electrode};
    \node at (0.5, 1-\hel/2) {top electrode};

  \fill[pattern=dots, pattern color=gray!50] (0,\hel) rectangle (1,1-\hel);
  \node at (0.5, 0.5) {$\Omega_{\mathrm{org}}$};

  \draw[thick, blue] (0,\hel) -- (1,\hel) node[midway, below, yshift=6mm] {$\Gamma_{\mathrm{bot/org}}$};
  \draw[thick, blue] (0,1-\hel) -- (1,1-\hel) node[midway, above, yshift=-6mm] {$\Gamma_{\mathrm{top/org}}$};
  \draw[thick, green!70!black] (0,0+\hel) -- (0,1-\hel) node[midway, left, xshift=-2.4mm, rotate=90] {$\Gamma_{\mathrm{ins}}$};
  \draw[thick, green!70!black] (1,0+\hel) -- (1,1-\hel) node[midway, right, xshift=2.4mm, rotate=90] {$\Gamma_{\mathrm{ins}}$};
\end{tikzpicture}
    \caption{Sketch of the domain covering the organic active layer. We differentiated between insulated boundaries $\Gamma_{\mathrm{ins}}$ and the boundary part adjacent to the top or bottom electrodes $\Gamma_{\mathrm{bot/org}} \cup \Gamma_{\mathrm{top/org}}$.}
    \label{fig_sketchdomain}
  \end{minipage}
\end{figure}

\subsection{Morphology formation}
The phasefield model describing the nanomorphology within the active layer is based on our previous work presented in \cite{PCiloglu_2025}. All details on the morphology model can be found in this reference, and we will therefore only briefly discuss the morphology evolution. The model traces the evolution of donor and acceptor domains within the fluid blend of donor, acceptor, and solvent. Phase separation is driven by the minimization of a free energy function $f$. In this paper, we distinguish between two different approaches: (i) a ternary donor–acceptor–solvent system and (ii) an extended system that additionally includes the surrounding air, explicitly incorporating solvent evaporation. In the ternary system (i), the evolution of polymer, NFA, and solvent volume fractions are governed by Cahn--Hilliard equations. Evaporation of solvent is imposed by an outward flux of solvent and respective inward flux of polymer and NFA at the top part of the domain. In the extended system (ii), we explicitly account for evaporative effects by means of additional phase fields. We distinguish between liquid and vapor phases and let evaporation, too, be driven by the free energy functional. While this second model certainly has its benefits, with a physics-based description of the evaporation effects, it has some disadvantages, like a large discretized systems and the need to postprocess the dried film height and subsequent assignment of the active layer domain $\Omega_\mathrm{org}$. For this reason, we will make use of both models, in Sections~\ref{sec_example1} and \ref{sec_example2}, respectively. 
\subsection{Exciton equation}
Excitons are created by incident photons within the organic active layer. Once generated, an exciton diffuses through the active layer until it recombines at the end of its lifetime, emitting a photon. If, however, the exciton reaches a donor--acceptor interface within its lifetime, the exciton dissociates, i.e., the hole and electron separate due to the difference in energy levels between the acceptor and donor materials. We define the exciton density ${X : \Omega_{\mathrm{org}} \times [0, t_{\mathrm{max}}] \to \mathbb{R}_+}$ and assume that excitons are generated with a constant rate $G$ due to incident photons. After generation, they can either recombine with rate $\eta_r$ or dissociate with rate $\eta_{d}$ at the donor--acceptor interface, indicated by the phase field gradient $\lambda |\nabla \phi |$. These effects lead to the following diffusion equation for the evolution of $X$: 
\begin{align} \label{eqn:exc}
\frac{\partial X}{\partial t} &= d_X \Delta X +  G - \eta_{r} X - \eta_{d} \lambda |\nabla \phi|X\quad \text{ in }\Omega_{\mathrm{org}} \times (0,t_{\mathrm{max}}).
\end{align}
The diffusion coefficient $d_X = \frac{l_X^2}{\tau_X}$ is connected to the exciton diffusion length $l_X$ and exciton lifetime $\tau_X$.
\subsection{Drift--Diffusion equations: Free charge carriers}
The starting point for modeling charge transport are the continuity equations for charges generated by electron and hole densities ${n: \Omega_{\mathrm{org}} \times [0, t_{\mathrm{max}}] \to \mathbb{R}_+}$ and ${p: \Omega_{\mathrm{org}} \times [0, t_{\mathrm{max}}] \to \mathbb{R}_+}$, 
\begin{align}\label{eqn:conteq_n}
    q_e \frac{\partial n}{\partial t} - \nabla \cdot \mathbf{j}_n &= - q_e R(n,p) + q_e\eta_{d} \lambda |\nabla \phi|X &\text{in }  \Omega_{\mathrm{org}} \times (0,t_{\mathrm{max}}), \\ \label{eqn:conteq_p}
    q_e \frac{\partial p}{\partial t} +\nabla \cdot \mathbf{j}_p &= - q_e R(n,p) + q_e \eta_{d} \lambda |\nabla \phi|X &\text{in }  \Omega_{\mathrm{org}} \times (0,t_{\mathrm{max}}). 
\end{align}
Note that both equations contain the same sink and source terms. More specifically, $R(n,p)$ denotes the recombination rate of electrons and holes, while $\lambda |\nabla \phi|X$ models the creation of free electrons and holes due to exciton dissociation at the donor--acceptor interface. The behavior of the charge fluxes $\mathbf{j}_n$ and $\mathbf{j}_n$ is highly dependent on the distribution of the energy levels of the material. Specifically, it is important to observe the difference in energetic disorder between organic and inorganic semiconductors. In inorganic cells, electronic states in either conductor and valance band are located at almost constant energies, giving rise to well-defined energy levels of conductor and valance band, respectively. In organic materials, random alignment of the molecules leads to disordered distributed energy levels. As a result, the energetic distribution in inorganic materials is often described via Boltzmann or Fermi--Dirac statistics \cite{markowich1990semiconductor}. Organic materials, on the other hand, require the use of Gaussian distributed energy levels, compare, e.g., \cite{burke2015device, Doan2019}. The number of transport states per energy level $E$ can hence be described by
\begin{align*}
    N_{\mathrm{\mathrm{Gauss}}}(E) = \frac{N_0}{\sigma \sqrt{2 \pi}} \exp\left( - \left( \frac{E - E_0}{\sqrt{2}\sigma} \right)^2 \right)
\end{align*}
with mean energy level $E_0$, variance $\sigma$, and total density of transport states $N_0$. The variance $\sigma$ is a material property, and is directly coupled to the degree of disorder within the organic material. The average energy level $E_0$ is, depending on the charge considered, coupled to the HOMO $E_H(\phi)$ or LUMO level $E_L(\phi)$ of the materials in addition to the contribution $q_e \psi$ of any external electric field with potential $\psi$. The Fermi function 
\begin{align*}
    f(E) = \exp\left( \frac{E-E_F}{k_B T} +1 \right)^{-1}
\end{align*}
assigns the probability of a charge being in a state with energy $E$ for a given Fermi level $E_F$. A key realization is that electrons and holes relax sufficiently fast on their respective levels such that their distributions can be considered independently. Defining the quasi-Fermi levels $\varphi_n$ and $\varphi_p$, for electrons and holes, respectively, their densities are given by 
\begin{align*}
    n &= \int_{-\infty}^{\infty} N_{\mathrm{\mathrm{Gauss}}, n}(E) f_n(E)  d E, \quad
    p = \int_{-\infty}^{\infty} N_{\mathrm{\mathrm{Gauss}}, p}(E) f_p(E)  d E,   
\end{align*}
with $N_{\mathrm{Gauss}, (n,p)}(E) f_{n,p}(E)$ the Gaussian distribution and Fermi function with adjusted mean and variance. Following the work presented in \cite{Doan2019, GPaasch_SScheinert_2010}, these equations can be rewritten to 
\begin{align} \label{eq_n}
n = N_{n0} \mathcal{G}\left( \frac{q_e \left(\psi-\varphi_n \right) - E_L}{k_B T}, \frac{\sigma_n}{k_B T} \right)  \quad \text{in} \ \Omega_{\mathrm{org}}\times (0,t_{\mathrm{max}}), \\ \label{eq_p}  p=N_{p0} \mathcal{G}\left( \frac{E_H - q_e \left( \psi - \varphi_p \right)}{k_B T}, \frac{\sigma_p}{k_B T} \right) \quad \text{in} \ \Omega_{\mathrm{org}}\times (0,t_{\mathrm{max}}), 
\end{align}
where the function $\mathcal{G}$ is given by the Gauss-Fermi integral, 
\begin{align} \label{eq_Gaussint}
\mathcal{G}(\eta, z) := \frac{1}{\sqrt{2\pi}} \int_{-\infty}^{\infty} \exp\left(-\frac{\xi^2}{2}\right) \frac{1}{\exp(z\xi-\eta)+1} d\xi.
\end{align}
At this point, we already note the expected difficulties with equations of the form \eqref{eq_Gaussint}, as is also extensively discussed in \cite{Doan2019, GPaasch_SScheinert_2010}. Therefore, we approximate the Gauss-Fermi integral by the Boltzmann approximation, compare \cite{GPaasch_SScheinert_2010}, valid in the low-density limit for small variances, 
\begin{align}\label{eq_GFapprox}
    \mathcal{G}(\eta, z) \approx \exp(\eta) \exp \left(\frac{z^2}{2}\right).
\end{align}
Consequently equations \eqref{eq_n}, \eqref{eq_p} reduce to
\begin{align}\label{eqn:elechole}
n &= N_{n0} \exp\left(\frac{\sigma_n^2}{2k_B^2 T^2}\right) \exp\left(\frac{q_e \left( \psi-\varphi_n \right)- E_L}{k_B T} \right)  \quad \text{in} \ \Omega_{\mathrm{org}}\times (0,t_{\mathrm{max}}), \\
p & = N_{p0} \exp\left(\frac{\sigma_p^2}{2k_B^2 T^2}\right) \exp\left(\frac{q_e \left(\varphi_p-\psi \right) + E_H}{k_B T} \right) \quad \text{in} \ \Omega_{\mathrm{org}}\times (0,t_{\mathrm{max}}). 
\end{align}
We note that this is a gross approximation, greatly simplifying the charge hopping transport typical for organic materials. However, to obtain feasible numerical methods, this approximation is deemed necessary. It remains to define the recombination rate of electrons and holes on the RHS of equation~\eqref{eqn:conteq_n}, \eqref{eqn:conteq_p}. We use the Langevin recombination rate $R(n,p) = \gamma \left(np- N^2_{intr}\right)$, valid for low mobility materials, see \cite{Deibel_2010}, 
with the constant Langevin recombination prefactor $\gamma$. The intrinsic carrier concentration $N_{intr}$ is given by 
\begin{align*}
    N_{intr}^2 = N_{n0}N_{p0} \exp \left( - \frac{E_L - E_H}{k_B T}\right). 
\end{align*}
The charge fluxes in equations~\eqref{eqn:conteq_n}, \eqref{eqn:conteq_p} are defined by
\begin{align}\label{eqn:chargeflux}
\mathbf{j}_n &= - \mu_n q_e n \nabla \varphi_n \quad \text{in }  \Omega_{\mathrm{org}} \times (0,t_{\mathrm{max}}), \\
\mathbf{j}_p &= - \mu_p q_e p \nabla \varphi_p \quad \  \text{in }  \Omega_{\mathrm{org}} \times (0,t_{\mathrm{max}}),
\end{align}
where $\mu_n$ and $\mu_p$ are the constant carrier mobilities. Note that, using equation~\eqref{eqn:elechole},
\begin{align*}
    -\nabla n + n \left( \nabla \frac{q_e \psi}{k_B T} - \nabla \frac{E_L}{k_B T} \right) = n \nabla \frac{q_e \varphi_n}{k_B T}, \\
    \nabla p + p \left( \nabla \frac{q_e \psi}{k_B T} - \nabla \frac{E_H}{k_B T} \right) = p \nabla \frac{q_e \varphi_p}{k_B T}, 
\end{align*}
and consequently, we can rewrite the fluxes to
\begin{align*}
    \mathbf{j}_n = k_B T \mu_n  \nabla n - \mu_n n \left(q_e \nabla \psi  - \nabla E_L \right) \hbox{ and } \ \mathbf{j}_p = - k_B T  \mu_p \nabla p - \mu_p p \left(q_e \nabla \psi- \nabla E_H \right).
\end{align*}
\begin{remark}
The present model explicitly resolves configurational disorder through bulk heterojunction morphologies, while neglecting energetic disorder associated with charge trapping. In particular, trapping and detrapping processes, commonly represented through Shockley–Read–Hall–type trap-state models, are not included. This simplification is introduced in order to make fully coupled drift–diffusion–exciton simulations feasible in a three-dimensional space and to provide a robust framework for exploring how configurational disorder affects charge transport and overall device performance. Explicit inclusion of trap occupancies would increase numerical stiffness, rendering the problem extremely difficult to solve and thus referred to future work.

We note that it is well established that energetic disorder plays a central role in organic photovoltaic devices, and Gaussian \cite{HBassler_1993} density-of-states models are often employed to describe the overall structure of the organic band manifold. In the present work, this Gaussian energetic landscape is effectively collapsed into a Maxwell--Boltzmann statistical description in order to obtain a tractable transport
formulation. However, deep tail states involved in recombination are frequently described using exponential \cite{JNelson_2003} distributions in order to reproduce ideality factors and transient charge-carrier dynamics. Trap states, therefore, remain essential in device models aimed at quantitatively reproducing experimental current–voltage characteristics, time-domain measurements, and recombination mechanisms.
\end{remark}
\subsection{Poisson equation: Electrical potential}
The electrical potential $\psi$ is coupled to the total free charge carrier density ($n+p)$ by
\begin{align}\label{eqn:Poisson}
- \nabla \cdot ( \varepsilon\nabla \psi) &= q_e \left( -n + p\right) \quad \text{in } \Omega_{\mathrm{org}}\times (0,t_{\mathrm{max}}),
\end{align}
where $\varepsilon = \varepsilon_0 \varepsilon_r$ with $\varepsilon_0$ and $\varepsilon_r$ are the dielectric constant and relative permittivity, respectively.
\subsection{Initial and boundary conditions} 
The governing equations are supplemented with initial conditions for the exciton density $X$, the electric potential $\psi$, and the quasi-Fermi levels $\varphi_n, \varphi_p$, 
\begin{align*}
    X(x, 0) = X^0(x),  \ \psi(x,0) = \psi^0(x) , \ \varphi_n (x,0) = \varphi_n^0(x) , \ \varphi_p(x,0) = \varphi_p^0 (x) \ \text{in } \Omega_{\mathrm{org}} .
\end{align*}
For the exciton density, we consider no-flux boundary conditions on the insulated part of the boundary as well as the boundary adjacent to the electrodes, 
\begin{align*}
\nabla X \cdot \bn &= 0 \quad \text{on }\partial\Omega_{\mathrm{org}} \times (0,t_\mathrm{max}) .
\end{align*}
For the electrostatic potential and the quasi-Fermi levels, we differentiate between the insulated boundary and the electron boundaries. On the insulated boundary, we impose no-flux conditions,
\begin{align*}
\varepsilon \nabla \psi \cdot \bn &= 0 \ & \text{on }&\Gamma_{\mathrm{ins}} \times (0,t_\mathrm{max}), \ &\mathbf{j}_n \cdot \bn & = 0 \quad & \text{on }&\Gamma_{\mathrm{ins}} \times (0,t_\mathrm{max}), \\
&&&& \mathbf{j}_p \cdot \bn &= 0 \quad & \text{on }&\Gamma_{\mathrm{ins}} \times (0,t_\mathrm{max}). 
\end{align*}
On the top and bottom boundary, we define the electrical potential satisfying the local electroneutrality condition $\psi^*$, 
\begin{align*}
    N_{n0} \mathcal{G}\left( \frac{q_e \psi^* - E_L}{k_B T}, \frac{\sigma_n}{k_B T} \right) = N_{p0} \mathcal{G}\left( \frac{E_H - q_e\psi^*}{k_B T}, \frac{\sigma_p}{k_B T} \right) \quad \text{on} \ \Gamma_{\mathrm{bot}} \cup \Gamma_{\mathrm{top}} \times (0,t_\mathrm{max}). 
\end{align*}
Equivalently, using the simplification~\eqref{eq_GFapprox}, this rewrites to
\begin{align*}
    &N_{n0} \exp \left(\frac{\sigma_n^2}{2k_B^2 T^2} \right) \exp\left(\frac{q_e \psi^* - E_L}{k_B T} \right) \\ & \qquad = N_{p0} \exp\left(\frac{\sigma_p^2}{2k_B^2 T^2}\right) \exp\left(\frac{E_H - q_e\psi^*}{k_B T} \right) \quad \text{on} \ \Gamma_{\mathrm{bot}} \cup \Gamma_{\mathrm{top}} \times (0,t_\mathrm{max}). 
\end{align*}
The boundary conditions for the electrostatic potential and quasi-Fermi levels are then given by, compare \cite{Glitzky2024stack, Doan2019}, 
\begin{align}\label{eq_bc_current}
\varphi_n &= V_{\mathrm{bot}} &&\text{on }\Gamma_{\mathrm{bot/org}}\times (0,t_\mathrm{max}),            \ &\varphi_n &= V_{\mathrm{top}} &&\text{on }\Gamma_{\mathrm{top/org}} \times (0,t_\mathrm{max}),\\
\varphi_p &= V_{\mathrm{bot}} &&\text{on }\Gamma_{\mathrm{bot/org}} \times (0,t_\mathrm{max}),           \ &\varphi_p &= V_{\mathrm{top}} &&\text{on }\Gamma_{\mathrm{top/org}} \times (0,t_\mathrm{max}),\nonumber \\ 
\psi &= V_{\mathrm{bot}} + \psi^{*} &&\text{on }\Gamma_{\mathrm{bot/org}}\times (0,t_\mathrm{max}),     \ &\psi &= V_{\mathrm{top}} + \psi^{*} &&\text{on }\Gamma_{\mathrm{top/org}} \times (0,t_\mathrm{max}),  \nonumber
\end{align}
where $V_{\mathrm{bot}}$ and $V_{\mathrm{top}}$ define the externally applied voltage difference.

\subsection{Device performance}
The electrical performance of organic devices is captured by current-density--voltage characteristic (I-V curves), where the device output current is mapped against the external voltage,
\begin{align} \label{eqn:current}
    I_{\theta}(V_{\mathrm{bot}}, V_{\mathrm{top}}) = \int_{\Gamma_{\theta}} (\mathbf{j}_n + \mathbf{j}_p) \cdot \bn \, d s, \quad \theta \in \{ \mathrm{top/org}, \mathrm{bot/org}\}\ .
\end{align}
From this characteristic, device specifics like the open-circuit voltage ($V_{OC}$), short-circuit current density ($J_{SC}$), and fill factor ($FF$) can be extracted. 
\section{Nondimensionalized governing equations} \label{sec:nondim}
We fix the typical length $x_c$, time $t_c$, temperature $T_c$, charge $q_c$, and voltage $V_c$ scales. Typical, dimensional quantities are denotes by $(\cdot)_c$, non-dimensional quantities by $\widehat{(\cdot)}$. Naturally, the typical charge is chosen to be the elementary charge $q_c = q_e$ and the typical temperature set to the constant temperature $T_c = T$. The typical voltage is defined to be the thermal voltage $V_c = \frac{k_B T_c}{q_e}$. Defining the nondimensional electron and hole densities $\widehat{X} = x_c^3 X$, $\widehat{n} = x_c^3 n$, $\widehat{p} = x_c^3 p$, the nondimensional charge fluxes $\widehat{\mathbf{j}}_n = \frac{x_c^2 t_c}{q_e} \mathbf{j}_n $, $\widehat{\mathbf{j}}_p = \frac{x_c^2 t_c}{q_e} \mathbf{j}_p$, and the nondimensional potentials $ \widehat{\psi} = \frac{\psi}{V_c}$, $\widehat{\varphi}_n = \frac{\varphi_n}{V_c}$, $\widehat{\varphi}_p = \frac{\varphi_p}{V_c}$, the nondimensional governing equations $\text{in }\widehat{\Omega}_{\mathrm{org}} \times (0,\widehat{t}_\mathrm{max})$ read
\begin{subequations}\label{eqn:gover}
\begin{align}
\frac{\partial \widehat{X}}{\partial \widehat{t}} &= \widehat{d}_X \widehat{\Delta} \widehat{X} +  \widehat{G} - \widehat{\eta_{r}} \widehat{X} -  \widehat{\eta}_{d} |{ \widehat{\nabla}} \widehat{\phi}|\widehat{X}, \label{eqn:nonExc} \\
\frac{\partial \widehat{n}}{\partial \widehat{t}} - \widehat{\nabla} \cdot \widehat{\mathbf{j}}_n  &= - \widehat{R} + \widehat{\eta}_d | \widehat{\nabla} \widehat{\phi}| \widehat{X}, \quad \widehat{\mathbf{j}}_n = \widehat{\mu}_n \widehat{\nabla} \widehat{n} - \widehat{\mu}_n \widehat{n} \left(\widehat{\nabla} \widehat{\psi} - \widehat{\nabla} \widehat{E}_L \right) , \label{eqn:nonElec} \\
\frac{\partial \widehat{p}}{\partial \widehat{t}} +\widehat{\nabla} \cdot \widehat{\mathbf{j}}_p &= -\widehat{R} + \widehat{\eta}_d | \widehat{\nabla} \widehat{\phi}| \widehat{X}, \quad \widehat{\mathbf{j}}_p = -\widehat{\mu}_p \widehat{\nabla} \widehat{p} - \widehat{\mu}_p \widehat{p} \left(\widehat{\nabla} \widehat{\psi} - \widehat{\nabla} \widehat{E}_H \right) , \label{eqn:nonHole} \\
\widehat{n} &= \widehat{N}_{n0} \exp\left(\frac{\widehat{\sigma}_n^2}{2}\right) \exp\left(\left( \widehat{\psi}-\widehat{\varphi}_n \right)- \widehat{E}_L\right) ,  \label{eqn:nonElecFer} \\
\widehat{p} & = \widehat{N}_{p0} \exp\left(\frac{\widehat{\sigma}_p^2}{2}\right) \exp\left(\left(\widehat{\varphi}_p-\widehat{\psi} \right) + \widehat{E}_H \right) , \label{eqn:nonHoleFer} \\
\widehat{R} &= \widehat{\gamma} \left(\widehat{n}\widehat{p}- \widehat{N}_{intr}^2 \right),  \\
- \widehat{\nabla} \cdot \widehat{\varepsilon} \widehat{\nabla} \widehat{\psi} &= - \widehat{n} + \widehat{p} . \label{eqn:nonPoisson}
\end{align}
\end{subequations}
The values of the nondimensional parameters are
\begin{align*}
    &\widehat{d}_X = \frac{t_c}{x_c^2}d_x,                  && \  \widehat{G} = t_c x_c^3 G,                                && \  \widehat{\eta}_{r} = t_c \eta_r,                           && \  \widehat{\eta}_{d} = \frac{t_c}{x_c} \lambda \eta_d,  \\
    &\widehat{\phi} (\widehat{x}) = \phi(x_c \widehat{x}),  && \  \widehat{\mu}_n = \frac{k_B T_c t_c}{x_c^2 q_e} \mu_n,    && \  \widehat{\mu}_p = \frac{k_B T_c t_c}{x_c^2 q_e} \mu_p,     && \  \widehat{N}_{n0} = x_c^3 N_{n0},  \\
    &\widehat{N}_{p0} = x_c^3 N_{p0},                       && \  \widehat{N}_{intr}^2 = x_c^6 N_{intr}^2,    && \  \widehat{\gamma} = \frac{t_c}{x_c^3} \gamma ,              && \  \widehat{\varepsilon} = \frac{k_B T_c x_c}{q_c^2} \varepsilon ,\\
    \noalign{$\widehat{E}_L(\widehat{\phi}) = \frac{E_{L}^{p} + (E_{L}^{\nfa} - E_{L}^{p}) \phi(x_c \widehat{x})}{k_B T_c} , \qquad \widehat{E}_H(\widehat{\phi}) = \frac{E_{H}^{p} + (E_{H}^{\nfa}- E_{H}^{p}) \phi(x_c \widehat{x})}{k_B T_c}. $} 
\end{align*}
The governing equations are equipped with the nondimensional initial conditions
\begin{align*}
    \widehat{X}(\widehat{x}, 0) &= \widehat{X}^0(\widehat{x})  && \text{in } \widehat{\Omega}_{\mathrm{org}} ,  \quad &\widehat{\psi}(\widehat{x}, 0) &= \widehat{\psi}^0(\widehat{x})  && \text{in } \widehat{\Omega}_{\mathrm{org}} ,\\
    \widehat{\varphi}_n (\widehat{x}, 0) &= \widehat{\varphi}_n^0(\widehat{x})  && \text{in } \widehat{\Omega}_{\mathrm{org}} , \quad &\widehat{\varphi}_p(\widehat{x}, 0) &= \widehat{\varphi}_p^0(\widehat{x})  && \text{in } \widehat{\Omega}_{\mathrm{org}} ,
\end{align*}
Neumann boundary conditions 
\begin{align*}
\widehat{\nabla} \widehat{X} \cdot \widehat{\bn} &= 0  && \text{on }\partial\widehat{\Omega}_{\mathrm{org}} \times (0,\widehat{t}_\mathrm{max}), \quad & \widehat{\mathbf{j}}_n \cdot \widehat{\bn} & = 0  &&\text{on }\widehat{\Gamma}_{\mathrm{ins}} \times (0,\widehat{t}_\mathrm{max}), \\ 
\widehat{\varepsilon} \widehat{\nabla} \widehat{\psi} \cdot \widehat{\bn} &= 0  && \text{on }\widehat{\Gamma}_{\mathrm{ins}} \times (0,\widehat{t}_\mathrm{max}), \quad  & \widehat{\mathbf{j}}_p \cdot \widehat{\bn} &= 0 && \text{on }\widehat{\Gamma}_{\mathrm{ins}} \times (0,\widehat{t}_\mathrm{max}),  
\end{align*}
and Dirichlet boundary conditions
\begin{align*}
\widehat{\varphi}_n &= \widehat{V}_\mathrm{bot} &&\text{on } \widehat{\Gamma}_{\mathrm{bot/org}}\times (0,\widehat{t}_\mathrm{max}),            \ &\widehat{\varphi}_n &= \widehat{V}_\mathrm{bot} &&\text{on } \widehat{\Gamma}_{\mathrm{top/org}} \times (0,\widehat{t}_\mathrm{max}),\\
\widehat{\varphi}_p &= \widehat{V}_\mathrm{bot} &&\text{on } \widehat{\Gamma}_{\mathrm{bot/org}} \times (0,\widehat{t}_\mathrm{max}),           \ &\widehat{\varphi}_p &= \widehat{V}_\mathrm{top} &&\text{on } \widehat{\Gamma}_{\mathrm{top/org}} \times (0,\widehat{t}_\mathrm{max}),\nonumber \\ 
\widehat{\psi} &= \widehat{V}_\mathrm{bot}+ \widehat{\psi}^{*} &&\text{on } \widehat{\Gamma}_{\mathrm{bot/org}} \times (0,\widehat{t}_\mathrm{max}),     \ &\widehat{\psi}&= \widehat{V}_\mathrm{top} + \widehat{\psi}^{*} &&\text{on } \widehat{\Gamma}_{\mathrm{top/org}} \times (0,\widehat{t}_\mathrm{max}),  \nonumber
\label{eqn:non_bc_current}
\end{align*}
where 
\begin{align*}
    & \widehat{N}_{n0} \exp \left(\frac{\widehat{\sigma}_n^2}{2} \right) \exp\left(\widehat{\psi}^* - \widehat{E}_L \right) \\
    & \qquad = \widehat{N}_{p0} \exp\left(\frac{\widehat{\sigma}_p^2}{2}\right) \exp\left(\widehat{E}_H - \widehat{\psi}^* \right)  \text{on} \ \widehat{\Gamma}_{\mathrm{bot}} \cup \widehat{\Gamma}_{\mathrm{top}} \times (0,\widehat{t}_\mathrm{max}). 
\end{align*}
The results presented in the remainder of this paper are based on the nondimensional equations, for readability, the hat-notation is dropped from now on.

\section{Discretization scheme}
In the following, the finite element method for approximating the solution of the governing equation given in Section \ref{sec:nondim} will be discussed. We will compare three different methods to solve nonlinear systems: Newton, Gummel, and Semi--Newton\-–Gummel. Although the Newton and Gummel schemes are commonly used in the literature, we develop a hybrid approach that combines the two methods to take advantage of their strengths.

\subsection{Spatial and temporal discretization}
Let $\{\mathcal{T}_h\}_h$ be a family of shape-regular triangulations of the organic domain $\Omega_{\mathrm{org}}$. Each mesh $\mathcal{T}_h$ consists of closed triangles such that ${\overline{\Omega_{\mathrm{org}}} = \bigcup_{K \in \mathcal{T}_h } \overline{K}}$ holds. For a given $\mathcal{T}_h$, we define the discretized function spaces 
\[
\mathcal{V}_h = \sett{ v \in C(\overline{\Omega_{\mathrm{org}}})}{ v \mid_{K}\in \mathbb{P}^\ell(K) \quad \forall K \in \mathcal{T}_h}  \subset H^1(\Omega_{\mathrm{org}})
\]
with $\mathbb{P}^\ell(K)$ being the set of all polynomials on $K$ of degree at most $\ell$. We seek approximations of the electrostatic potential $\psi_h \in \mathcal{V}^h$, the quasi-Fermi potentials $\varphi_{n,h}, \varphi_{p,h} \in \mathcal{V}^h$, and the exciton density $X_h \in \mathcal{V}^h$. The governing equations in the weak form read 
\begin{subequations}
\begin{equation}\label{eqn:discpossion}
\int_{\Omega_{\mathrm{org}}} \varepsilon \nabla \psi_h \cdot \nabla v \, dx = \int_{\Omega_{\mathrm{org}}} (- n_h + p_h) v \, dx
\end{equation}
\begin{align}\label{eqn:discelec}
\int_{\Omega_{\mathrm{org}}}  \frac{\partial n_h}{\partial t} v_n \, dx + & \int_{\Omega_{\mathrm{org}}} \mu_n \nabla n_h \cdot \nabla v_n \, dx + \int_{\Omega_{\mathrm{org}}}  \mu_n n_h ( \nabla \psi_h - \nabla E_L) \cdot \nabla v_n \, dx \nonumber \\
   & = -\int_{\Omega_{\mathrm{org}}} R_h v_n \, dx + \int_{\Omega_{\mathrm{org}}} \eta_{d}|\nabla \phi|X_h v_n \, dx
\end{align}
\begin{align}\label{eqn:dischole}
\int_{\Omega_{\mathrm{org}}}  \frac{\partial p_h}{\partial t} v_p \, dx + & \int_{\Omega_{\mathrm{org}}} \mu_p \nabla p_h \cdot \nabla v_p \, dx - \int_{\Omega_{\mathrm{org}}}  \mu_p p_h  (\nabla \psi - \nabla E_H) \cdot \nabla v_p\, dx \nonumber \\
   & = -\int_{\Omega_{\mathrm{org}}} R_h v_p \, dx + \int_{\Omega_{\mathrm{org}}} \eta_{d} |\nabla \phi|X_h v_p \, dx
\end{align}
\begin{align}\label{eqn:discexc}
\int_{\Omega_{\mathrm{org}}} \frac{\partial X_h}{\partial t} v_x \, dx&+ \int_{\Omega_{\mathrm{org}}} d_X \nabla X_h \cdot \nabla v_x \, dx \nonumber \\ 
&+ \int_{\Omega_{\mathrm{org}}} \left[(\eta_{r}+ \eta_{d}|\nabla \phi|) X_h - G\right] v_x \, dx = 0,
\end{align}
for the test functions $v, v_n, v_p, v_x \in \mathcal{V}^h $. Note that the carrier densities $n_h$ and $p_h$ are dependent on the Fermi-potentials $\varphi_{n,h}$, $\varphi_{p,h}$ and the electrical potential $\psi_h$. We define the discrete reaction term $R_h = R(n_h,p_h)$. Furthermore, we assume that the test functions $v, v_n, v_p, v_x$ vanish at the Dirichlet boundary of the domain and that the ansatz functions represented by $\psi_h, \varphi_{n,h}, \varphi_{p,h}, X_h$ satisfy the Dirichlet boundary conditions \eqref{eq_bc_current}.
\end{subequations}
For the time discretization we define an equidistant sequence of $0=t_0,\, t_1, \, \ldots \, ,$ $ t_n = t_\mathrm{max} $ with fixed time step size $\tau$. The fully discrete form of the governing equations \eqref{eqn:discpossion}--\eqref{eqn:discexc} is obtained using an implicit (backward) Euler scheme, whereafter the resulting nonlinear system is solved by the Newton method. The fully discrete form can be written as: 
\begin{subequations}
\begin{equation} \label{eqn:fullypossion}
\left( \varepsilon \nabla \psi_h^{k+1}, \nabla v\right)_{\Omega_{\mathrm{org}}} =  \left(- n_h^{k+1} + p_h^{k+1}, v\right)_{\Omega_{\mathrm{org}}},
\end{equation}
\begin{align}\label{eqn:fullyelec}
\left( \frac{n_h^{k+1} - n_h^{k}}{\tau}, v_n\right)_{\Omega_{\mathrm{org}}} &+  \left( \mu_n \nabla n_h^{k+1} , \nabla v_n\right)_{\Omega_{\mathrm{org}}} +\left( \mu_n n_h^{k+1} ( \nabla \psi_h^{k+1} - \nabla E_L),  \nabla v_n\right)_{\Omega_{\mathrm{org}}} \nonumber\\
& +\left(R_{h}^{k+1}  , v_n\right)_{\Omega_{\mathrm{org}}} =  \left( \eta_{d} |\nabla \phi|X_h^{k+1}   , v_n\right)_{\Omega_{\mathrm{org}}},
\end{align}
\begin{align}\label{eqn:fullyhole}
\left( \frac{p_h^{k+1} - p_h^{k}}{\tau}, v_p\right)_{\Omega_{\mathrm{org}}} &+  \left( \mu_p \nabla p_h^{k+1} , \nabla v_p\right)_{\Omega_{\mathrm{org}}} -\left( \mu_p p_h^{k+1} ( \nabla \psi_h^{k+1} - \nabla E_H),  \nabla v_p\right)_{\Omega_{\mathrm{org}}} \nonumber\\
& +\left(R_{h}^{k+1}  , v_n\right)_{\Omega_{\mathrm{org}}} =  \left( \eta_{d} |\nabla \phi|X_h^{k+1} , v_p\right)_{\Omega_{\mathrm{org}}},
\end{align}
\begin{align}\label{eqn:fullyexc}
\left( \frac{X_h^{k+1} - X_h^{k}}{\tau}, v_x\right)_{\Omega_{\mathrm{org}}}
+\left( d_X \nabla X_h^{k+1}, \nabla v_x\right)_{\Omega_{\mathrm{org}}}
 & +\left( (\eta_{r} + \eta_{d}|\nabla \phi|)  X_h^{k+1}, v_x\right)_{\Omega_{\mathrm{org}}} \nonumber \\
 & = \left( G, v_x\right)_{\Omega_{\mathrm{org}}}.
\end{align}
Here, $n_h^{k} = n_h^{k}(\varphi_{n,h}^{k}, \psi_h^{k+1})$, $p_h^{k} = p_h^{k}(\varphi_{p,h}^{k}, \psi_h^{k+1})$, $n_h^{k+1} = n_h^{k+1}(\varphi_{n,h}^{k+1}, \psi_h^{k+1})$, $p_h^{k+1} = p_h^{k+1}(\varphi_{p,h}^{k+1}, \psi_h^{k+1})$, and $R_{h}^{k+1} = R_{n,h}^{k+1}(n_h^{k+1},p_h^{k+1})$, allowing us to solve solely for the electrical potential $\psi_h^{k+1}$, the Fermi potentials $\varphi_{n,h}^{k+1}, \varphi_{p,h}^{k+1}$, and the exciton density $X_{h}^{k+1}$.
\end{subequations}
\subsection{Newton Method}
We discretize the system using an implicit scheme, which results in a single, fully-coupled, large, nonlinear system. This system is solved iteratively using Newton's method, which requires the solution of one large linear system at each iteration. We rewrite the Poisson equation \eqref{eqn:fullypossion} and the hole/electron equations \eqref{eqn:fullyelec}--\eqref{eqn:fullyhole} to its residual form
\begin{align}
\boldsymbol{F} = 
\begin{bmatrix}
    \boldsymbol{F}_{\psi}(\boldsymbol{\psi}^{k+1},\boldsymbol{\varphi_n}^{k+1},\boldsymbol{\varphi_p}^{k+1}) \\ 
    \boldsymbol{F}_{\varphi_n} (\boldsymbol{\psi}^{k+1},\boldsymbol{\varphi_n}^{k+1},\boldsymbol{\varphi_p}^{k+1}, \boldsymbol{X}^{k+1})  \\
    \boldsymbol{F}_{\varphi_p} (\boldsymbol{\psi}^{k+1},\boldsymbol{\varphi_n}^{k+1},\boldsymbol{\varphi_p}^{k+1}, \boldsymbol{X}^{k+1}) 
\end{bmatrix} = \boldsymbol{0}.
\end{align}
Here, the coefficient vectors $\boldsymbol{\psi}^{k}, \boldsymbol{\varphi_n}^{k}, \boldsymbol{\varphi_p}^{k}$ denote the discretized electrostatic potential and the quasi-Fermi levels at time $t_k$, respectively. Starting from a good initial guess, Newton's method constructs the following sequences 
\begin{align}
{\underbrace{\begin{bmatrix}
\frac{\partial \boldsymbol{F}_{\psi}}{\partial \boldsymbol{\psi}} & \frac{\partial \boldsymbol{F}_{\psi}}{\partial \boldsymbol{\varphi}_n} & \frac{\partial \boldsymbol{F}_{\psi}}{\partial \boldsymbol{\varphi}_p} \\
\frac{\partial \boldsymbol{F}_{\varphi_n}}{\partial \boldsymbol{\psi}} & \frac{\partial \boldsymbol{F}_{\varphi_n}}{\partial \boldsymbol{\varphi}_n} & \frac{\partial \boldsymbol{F}_{\varphi_n}}{\partial \boldsymbol{\varphi}_p} \\
\frac{\partial \boldsymbol{F}_{\varphi_p}}{\partial \boldsymbol{\psi}} & \frac{\partial \boldsymbol{F}_{\varphi_p}}{\partial \boldsymbol{\varphi}_n} & \frac{\partial \boldsymbol{F}_{\varphi_p}}{\partial \boldsymbol{\varphi}_p}
\end{bmatrix}}_{\mathbf{J}^{k+1,m}_{\boldsymbol{F}}} }^{k+1,m}
\begin{pmatrix}
\delta \boldsymbol{\psi}^{k+1,m} \\
\delta \boldsymbol{\varphi}_n^{k+1,m}\\
\delta \boldsymbol{\varphi}_p^{k+1,m}
\end{pmatrix}
= 
-
\begin{bmatrix}
\boldsymbol{F}_{\psi} \left( \boldsymbol{\psi}^{k+1,m}, \boldsymbol{\varphi}_n^{k+1,m}, \boldsymbol{\varphi}_p^{k+1,m} \right) \\
\boldsymbol{F}_{\varphi_n} \left( \boldsymbol{\psi}^{k+1,m}, \boldsymbol{\varphi}_n^{k+1,m}, \boldsymbol{\varphi}_p^{k+1,m} \right) \\
\boldsymbol{F}_{\varphi_p} \left( \boldsymbol{\psi}^{k+1,m}, \boldsymbol{\varphi}_n^{k+1,m}, \boldsymbol{\varphi}_p^{k+1,m} \right)
\end{bmatrix}
\end{align}
for the $(m+1)$th Newton update vector 
\begin{subequations}
\begin{align}
\boldsymbol{\psi}^{k+1,m+1} &= \boldsymbol{\psi}^{k+1,m} + \delta\boldsymbol{\psi}^{k+1,m},  \\
\boldsymbol{\varphi_n}^{k+1,m+1} &= \boldsymbol{\varphi_n}^{k+1,m} + \delta\boldsymbol{\varphi_n}^{k+1,m}, \\
\boldsymbol{\varphi_p}^{k+1,m+1} &= \boldsymbol{\varphi_p}^{k+1,m} + \delta\boldsymbol{\varphi_p}^{k+1,m}, 
\end{align}
\end{subequations}
where $\mathbf{J}^{k,m}_{\boldsymbol{F}}$ is the Newton derivative (Jacobian matrix) of $\boldsymbol{F}^{k,m}$. Lastly, the fully discrete exciton equation \eqref{eqn:fullyexc} can be expressed in the following form:
\begin{align}
    \frac{1}{\tau} \mathbf{M}\left(\boldsymbol{X}^{k+1} - \boldsymbol{X}^{k}\right) + d_X \mathbf{K} \boldsymbol{X}^{k+1} +  \left(\eta_{r} + \eta_{d}|\nabla \phi|\right)\mathbf{M}\boldsymbol{X}^{k+1} = \mathbf{f_X}.
\end{align}
By rewriting, we obtain the linear system
\begin{align}\label{eqn:linexc}
   \left( \mathbf{M} + d_X \mathbf{K} + \left(\eta_{r} + \eta_{d} |\nabla \phi|\right)\mathbf{M} \right)   \boldsymbol{X}^{k+1}  = \mathbf{f_X} +  \mathbf{M} \boldsymbol{X}^{k}, 
\end{align}
where the coefficient vector $\boldsymbol{X}^{k}$ represents the discretized exciton density. The matrices $\mathbf{M}$ and $\mathbf{K}$ denote the standard mass and stiffness matrices arising from the spatial discretization. The right-hand side vector $\boldsymbol{f_X}$ represents the discretized representations of terms containing $G$. 

The Jacobian matrix $\mathbf{J}_{\boldsymbol{F}}$ is non-symmetric due to the presence of the drift term. The Generalized Minimal Residual (\texttt{GMRES}) \cite{saad2003iterative} method is the solver of choice for non-symmetric linear systems. \texttt{GMRES} minimizes the residual norm over the Krylov subspace generated by the matrix $\mathbf{J}_{\boldsymbol{F}}$, providing predictable convergence for many such systems. Despite its robustness, the convergence rate of \texttt{GMRES} might drastically deteriorate for large, ill-conditioned matrices that are frequently used in the simulation of semiconductor devices. In order to accelerate the convergence and improve the conditioning of the system, we employ the Algebraic Multigrid (\texttt{AMG}) method as the primary preconditioner. The final step in the sequential solution process involves updating the exciton density $\boldsymbol{X}$ by solving the linear system \eqref{eqn:linexc}. As the system could be solved using a symmetric Krylov solver as well as a direct solver, we continue to adhere to the general numerical framework by utilizing the \texttt{GMRES} iterative solver with an \texttt{AMG} precondition.

\subsection{Gummel Method}
While Newton's method provides quadratic convergence, its requirement for a suitable initial solution and the need for forming and solving a large, non-symmetric Jacobian matrix can be prohibitively expensive. To enhance computational efficiency, the final non-linear system is instead resolved using the Gummel-type iteration \cite{Brinkman2013modell_asymptotics, Gummel_1964}, a semi-implicit approach that exploits the structural properties of the governing equations.

We first determine the electrical potential $\psi_h^{k+1}$ from the previous states of $n_{h}^{k}$, $p_{h}^{k}$, $X_{h}^{k}$, then we update $n_{h}^{k+1}, p_{h}^{k+1}, X_{h}^{k+1}$ by using the updated electrical potential $\psi_h^{k+1}$. The fully discrete form of the Poisson equation can be written as:
\begin{align}\label{eqn:linpoisson}
     \varepsilon \mathbf{K} \boldsymbol{\psi}^{k+1}  = \mathbf{f}_{\psi}^k. 
\end{align}
The fully discretized electron and hole equations are expressed as
\begin{subequations}
\begin{align} \label{eqn:linelecholn}
    \mu_n \mathbf{K} \boldsymbol{n}^{k+1} +  \mu_n \mathbf{C}_n\boldsymbol{n}^{k+1} + \mathbf{R}_n \boldsymbol{n}^{k+1} &= \mathbf{f}_n^k, \\ \label{eqn:linelecholp}
    \mu_n \mathbf{K} \boldsymbol{p}^{k+1} +  \mu_n \mathbf{C}_p\boldsymbol{p}^{k+1} + \mathbf{R}_p \boldsymbol{p}^{k+1} &= \mathbf{f}_p^k, 
\end{align}
where $\mathbf{C_n} = (\mathbf{c}_{ij}^n) $, $\mathbf{C_p} = (\mathbf{c}_{ij}^p) $, $\mathbf{R_n} = (\mathbf{r}_{ij}^n) $, and $\mathbf{R_p} = (\mathbf{r}_{ij}^p) $ with the individual terms given by 
\begin{align*}
    \mathbf{c}^n_{ij} &= \left( ( \nabla \psi_h^{k+1} - \nabla E_L)\zeta^i, \nabla \zeta^j \right), & \mathbf{r}^n_{ij} &=\left( \zeta^i n_h^k,  \zeta^j \right) \\
    \mathbf{c}^p_{ij} &= \left( ( \nabla \psi_h^{k+1} - \nabla E_H)\zeta^i, \nabla \zeta^j \right), & \mathbf{r}^p_{ij} &= \left( \zeta^i p_h^k , \zeta^j  \right)
\end{align*}
for $\{\zeta^i\} \subset \mathcal{V}^h$ . The right-hand side vector $\mathbf{f}_{\psi}^k$ is the discretized representation of the terms involving $n_h^k$ and $p_h^k$, while the vectors $\mathbf{f}_n^k$, and $\mathbf{f}_p^k$ represent the discretized form of terms $ \eta_{d} |\nabla \phi| X_h^k $. The matrix form of the exciton equation is discussed in the previous section. 
\end{subequations}

To enhance the robustness of the Gummel iteration at high applied potentials, a damping strategy is introduced to prevent divergence caused by strong coupling between the electric potential and carrier densities; see, e.g., \cite{Brinkman2013modell_asymptotics}. In this approach, a damping parameter $\alpha \in (0,1)$ is employed to form a convex combination of the current and previous iterates, thereby stabilizing the update for all variables. For each quantity $ u \in \{\psi, n, p, X\} $, the damped iteration step is given by
\begin{align}
    u^{k+1} = \alpha  u' + (1-\alpha)u^{k},
\end{align}
where $u'$ is the solution of each of the equations \eqref{eqn:linexc}, \eqref{eqn:linpoisson}, \eqref{eqn:linelecholn}, and \eqref{eqn:linelecholp}. Note that the boundary conditions are preserved, since any convex combination of two admissible solutions remains admissible. The time iterative procedure is carried out until convergence is reached, as measured by the $L^2$-norm of the difference between consecutive iterations. After reaching the stopping criteria, the external voltage difference $ V_{\mathrm{top}} - V_{\mathrm{bot}}$ is slightly increased, and the iteration is continued until convergence is achieved. The whole procedure is summarized in the Algorithm \ref{alg:iter}. Since the Poisson equation is symmetric but the coupled electron-hole transport equations are not, the \texttt{GMRES} algorithm is chosen as the Krylov iterative solver to ensure consistency in the numerical treatment of all subsystems. An \texttt{AMG} preconditioner is incorporated to improve convergence.
\begin{algorithm}
\caption{Iteration with Voltage Update}
\small
\begin{algorithmic}[1]\label{alg:iter}
\FOR{$V_\mathrm{top},\, V_\mathrm{bot}$ up to final values}
    \STATE Set boundary voltages $V_\mathrm{top}$, $V_\mathrm{bot}$
    \STATE Initialize $\psi^0$, $\varphi_n^0$, $\varphi_p^0$, $X^0$
    \STATE Compute initial $n^0$ and $p^0$ from \eqref{eqn:nonElecFer} - \eqref{eqn:nonHoleFer}
\WHILE{$L^2$ error $>$ \texttt{tol}}
    \STATE Solve \eqref{eqn:fullypossion} for the electrical potential $\psi'$
    \STATE Update $\psi^{k+1}$ by the damping strategy with $\psi'$ 
    \STATE Solve \eqref{eqn:fullyelec}-\eqref{eqn:fullyhole} for electron and hole density $n',p' $ using current $\psi^{k+1}$, $n^{k}$, $p^{k}, X^k$
    \STATE Solve \eqref{eqn:fullyexc} for the exciton density $X^{k+1}$
    \STATE Update $n^{k+1},p^{k+1}, X^{k+1}$ by the damping strategy with $n', p', X'$ 
    \STATE Compute $L^2$ error between consecutive steps
\ENDWHILE
\STATE Increase $V_\mathrm{top}$ and $V_\mathrm{bottom}$
\ENDFOR
\end{algorithmic}
\end{algorithm}
\subsection{Semi--Newton--Gummel Method}
We use a third strategy that combines the advantages of both, Newton's method and the standard Gummel iteration. By using Newton's approach locally, the goal is to improve nonlinear convergence inside each subproblem while preserving the robustness of the Gummel decoupling. Because the nonlinear coupling between the electrostatic potential and the electron and hole quasi-Fermi levels often leads to convergence problems when the system is treated fully coupled, this hybrid approach is ideal for drift--diffusion systems.

We decouple the discretized system using the Gummel scheme and the standard FE approach, resulting in three nonlinear systems and one linear system. We rewrite the Poisson equation \eqref{eqn:fullypossion}, the electron/hole equations \eqref{eqn:fullyelec}--\eqref{eqn:fullyhole} in the residual form
\begin{subequations}
\begin{align}
\boldsymbol{F}_{\psi}(\boldsymbol{\psi}^{k+1}, \boldsymbol{\varphi_n}^{k}, \boldsymbol{\varphi_p}^{k}) &= 0,  \\
\boldsymbol{F}_{\varphi_n}(\boldsymbol{\psi}^{k}, {\boldsymbol\varphi_n}^{k+1}, {\boldsymbol\varphi_p}^{k}, X^{k}) &= 0,  \\
\boldsymbol{F}_{\varphi_p}({\boldsymbol\psi}^{k}, {\boldsymbol\varphi_n}^{k}, {\boldsymbol\varphi_p}^{k+1}, X^{k}) &= 0. 
\end{align}
\end{subequations}
Instead of solving each subproblem with a simple fixed-point update, we apply a Newton iteration to each equation separately. The following sequences are constructed by Newton's method 
\begin{align*}
\mathbf{J}_{\boldsymbol{F}_{\psi}}({\boldsymbol\psi}^{k+1,m}, \boldsymbol{\varphi_n}^{k}, \boldsymbol{\varphi_p}^{k}) \, \delta\boldsymbol{\psi}^{k+1}
&= - \boldsymbol{F}_{\psi}(\boldsymbol{\psi}^{k+1,m}, \boldsymbol{\varphi_n}^{k}, \boldsymbol{\varphi_p}^{k}),  \\
\boldsymbol{\psi}^{k+1,m+1} &= \boldsymbol{\psi}^{k+1,m} + \delta\boldsymbol{\psi}^{k+1}, \qquad m = 1,\dots 
\end{align*}
\begin{align*}
\mathbf{J}_{\boldsymbol{F}_{\varphi_n}}(\boldsymbol{\psi}^{k+1}, \boldsymbol{\varphi_n}^{k+1,m}, \boldsymbol{\varphi_p}^{k}) \, \delta\boldsymbol{\varphi_n}^{k+1}
&= - \boldsymbol{F}_{\varphi_n}(\boldsymbol{\psi}^{k+1}, \boldsymbol{\varphi_n}^{k+1,m}, \boldsymbol{\varphi_p}^{k}),  \\
\boldsymbol{\varphi_n}^{k+1,m+1}
&= \boldsymbol{\varphi_n}^{k+1,m} + \delta\boldsymbol{\varphi_n}^{k+1}, \qquad m = 1,\dots 
\end{align*}
\begin{align*}
\mathbf{J}_{\boldsymbol{F}_{\varphi_p}}(\boldsymbol{\psi}^{k+1}, \boldsymbol{\varphi_n}^{k}, \boldsymbol{\varphi_p}^{k+1,m}) \, \delta\boldsymbol{\varphi_p}^{k+1}
&= - \boldsymbol{F}_{\varphi_p}(\boldsymbol{\psi}^{k+1}, \boldsymbol{\varphi_n}^{k}, \boldsymbol{\varphi_p}^{k+1,m}),\\
\boldsymbol{\varphi_p}^{k+1,m+1}
&= \boldsymbol{\varphi_p}^{k+1,m} + \delta\boldsymbol{\varphi_p}^{k+1}, \qquad m = 1,\dots 
\end{align*}
Similar to the Gummel iteration, we employed a damping strategy for the updates $\psi$, $\varphi_n$, $\varphi_p$, and $X$ after solving each equation. 

In analogy with the fully coupled Newton method and the classical Gummel iteration, each linear system arising in the decoupled Newton update is solved using \texttt{GMRES} with an \texttt{AMG} preconditioner. This ensures that all three approaches utilize the same robust linear solver, allowing for a fair comparison between the different nonlinear strategies.

\subsection{Approximation of terminal current}
In the simulations of semiconductor devices, the relation between the total current and the applied voltage is the key quantity of interest. Accurately determining the terminal current at electrical contacts is thus essential for post-processing simulation data. In the following, we outline a method for evaluating estimated terminal currents utilizing the finite element schemes; see, e.g., \cite{PFarrell_DPeschka_2019, PFarrell_NRotundo_DHDoan_MKantner_2017, QZhang_QWang_LZhang_BLu_2022}. 
For a given ohmic contact $\theta$, i.e., the non-isolated contact to the electrodes $\theta \in \{\mathrm{top/org}, \mathrm{bot/org} \}$, we define an auxiliary function $w_\theta$ such that $w_\theta: \Omega_{\mathrm{org}} \rightarrow \mathbb{R}$ is a scalar function defined by the solution of the following boundary value problem:
\begin{align*}
- \nabla^2 w_\theta &= 0 \quad \text{in } \Omega_{\mathrm{org}},  \\
\nabla w_\theta \cdot \mathbf{n} &= 0 \quad \text{on } \Gamma_{\mathrm{ins}}, \\
w_\theta &= \delta_{\theta,D} \quad \text{on } \Gamma_{D}, 
\end{align*}
where $\delta_{\theta, D}$ is the Kronecker delta. To represent the outflow current, we first take the time derivative of the Poisson equation \eqref{eqn:nonPoisson}
\begin{align*}
- \nabla \cdot \varepsilon \nabla (\partial_t \psi) &= -\partial_t n +\partial_t p   = - \nabla \cdot \mathbf{j}_n + R -   \eta_{d}|\nabla \phi|X   - \nabla \cdot \mathbf{j}_p - R + \eta_{d}|\nabla \phi|X  \\
&= - \nabla \cdot \mathbf{j}_n  - \nabla \cdot \mathbf{j}_p 
\end{align*}
Using this auxiliary function, the terminal current \eqref{eqn:current} through contact $\theta$ becomes
\begin{align*}
I_\theta &=  \int_{\partial \Omega_{\mathrm{org}}} w_\theta (\mathbf{j}_n+\mathbf{j}_p)  \cdot \mathbf{n} \, ds = \int_{\Omega_{\mathrm{org}}} \nabla \cdot (w_\theta (\mathbf{j}_n+\mathbf{j}_p) ) \, dx, \\
        &= \int_{\Omega_{\mathrm{org}}} \nabla w_\theta \cdot (\mathbf{j}_n + \mathbf{j}_p) \, dx - \int_{\Omega_{\mathrm{org}}} \nabla w_\theta \cdot \varepsilon \nabla(\partial_t \psi) \, dx.
\end{align*}
Applying the semi-discrete time approximation, the total terminal current $I_\theta$ flowing out of $\Gamma_\theta$ at the $k+1$th time step is given by
\[
I_\theta^{k+1} = \int_{\Omega_{\mathrm{org}}} \nabla w_\theta \cdot (\mathbf{j}_n^{k+1} + \mathbf{j}_p^{k+1}) \, dx + \frac{1}{\tau} \int_{\Omega_{\mathrm{org}}} \nabla w_\theta \cdot \varepsilon ( \nabla\psi^{k+1} - \nabla\psi^{k}) \, dx.
\]

\section{Numerical results} 
In this section, we provide several numerical results for two-dimensional (2D) and three-dimensional (3D) bulk heterojunction  architectures. The values for the various dimensionless parameters follow \cite{Brinkman2013modell_asymptotics}, and summarized in Table~\ref{tab:para}. 

\begin{table}[htp!]
\centering
\caption{Values of parameters used in the simulations.}
\begin{tabular}{c c c c}
\hline
  $d_X = 10^{-2}$ & $\mu_n = 3 $ & $\mu_p = 1 $  & $\varepsilon = 10^{-1}$\\
  $G = 16990$ & $\gamma = 0.6987$ & $T= 300$ & $\tau = 10^{-4}$\\
  $N_{n0}= 1$ & $N_{p0} = 1$ & $\eta_r = 1$ & $\eta_d = 1$ \\
  $\sigma_n = 0.3868$ & $\sigma_p = 0.3868$  & $N_{intr} = 0 $ & $\texttt{tol} = 10^{-4}$\\
  $E_L^{\nfa} = -4.10 $ & $E_L^p = -3.28$ & $E_H^{\nfa} = -5.65 $ & $E_H^p = -5.13$ \\
\hline
\end{tabular}
\label{tab:para}
\end{table}

\subsection{Example 1:} \label{sec_example1}
In this section, we present numerical solutions of the governing equations of \eqref{eqn:gover} using morphologies generated by the three-component phase-field model in 2- and 3-dimensional domains. The morphologies are generated on meshes with grid points, $n_x$, $n_y$, $n_z$. These uniform meshes in two and three dimensions of nondimensional size $10 \times 10$, and $10 \times 10 \times 10$, respectively, are generated with $n_x \times n_y$ or $n_x \times n_y \times n_z$ nodes, using linear Lagrange elements on triangular or tetrahedral cells. The same mesh discretizations are subsequently employed for the drift--diffusion simulations. The morphologies were generated using the same parameter set reported in \cite{PCiloglu_2025}. Unless stated otherwise, we set the Newton and \texttt{GMRES} solver tolerances to $\texttt{rtol} = 10^{-6}$ and $\texttt{atol} = 10^{-10}$, while the \texttt{AMG} preconditioner tolerance is set to $10^{-6}$. 

\begin{figure}[htp!]
    \centering
    \includegraphics[width=1\linewidth]{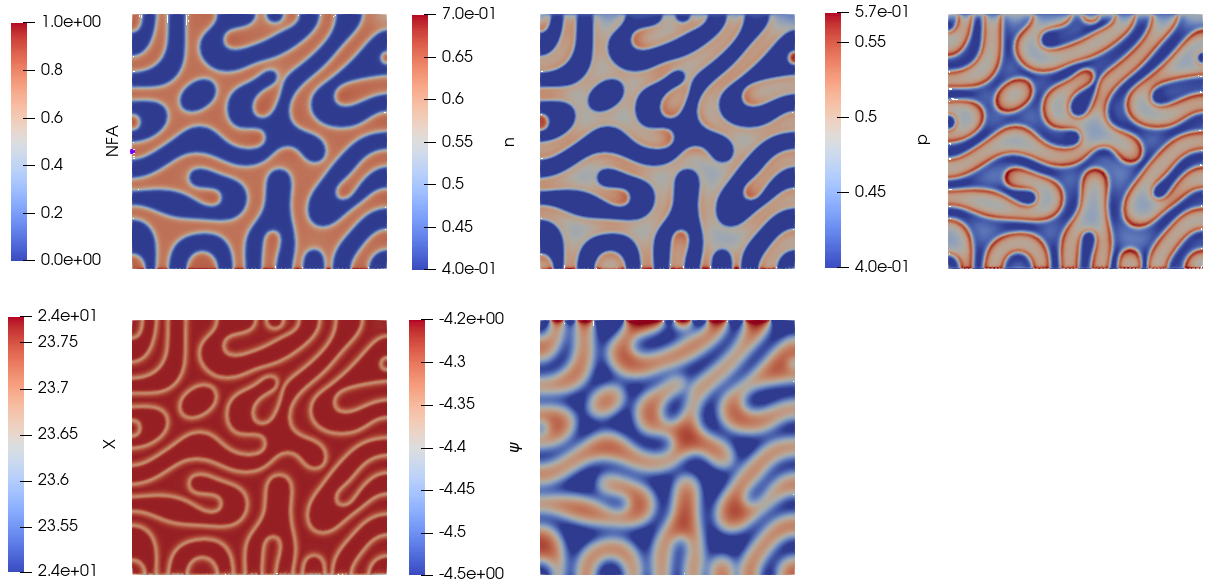}
    \caption{Morphology (left) and snapshots of the electron density (middle) and hole density (right) at time $t = 0.001$ in 2D, computed using the Newton method.}
    \label{fig:Ex1_Newtonnp}
\end{figure}

We use $n_x = 200$ and $n_y = 200$ grid points in 2D, and $n_x = 25$, $n_y = 25$, and $n_z = 25$ grid points in 3D. Our iteration starts with the external potentials $V_{\mathrm{top}} = 0$ and $V_{\mathrm{bot}} = 0$. After convergence, we increase $V_{\mathrm{top}}$ by $0.02$ while keeping $V_{\mathrm{bot}}$ fixed at $0$, and we repeat this procedure until $V_{\mathrm{top}} = 10$ is reached.
The model is solved with a \texttt{PYTHON} implementation using the finite element library \texttt{DOLFINX} \cite{dolfinx} from the \texttt{FENICS} project \cite{fenics} with version 0.9.0.

\begin{figure}[htp!]
    \centering
    \includegraphics[width=0.9\linewidth]{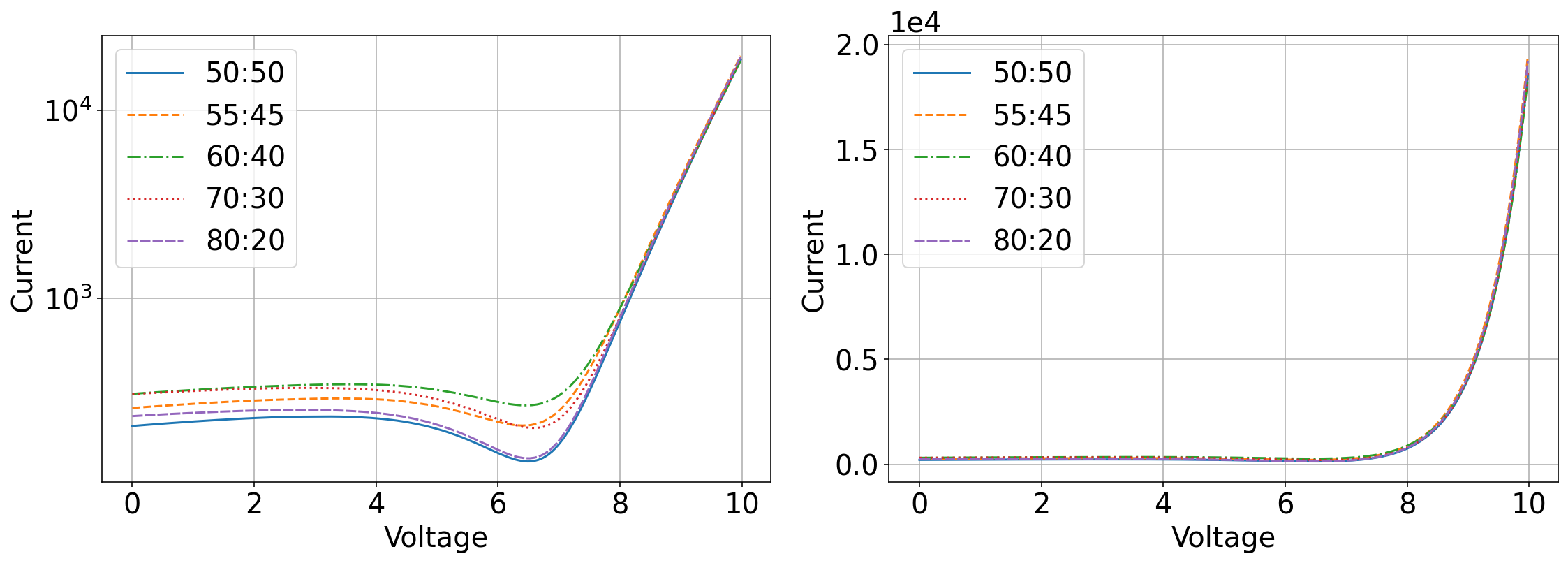}
    \caption{Total currents flowing out of the boundary $\Gamma_{\mathrm{bot/org}}$: semilogarithmic plot (left) and Cartesian plot (right), computed using 2D morphologies with different blend ratios $50{:}50$, $55{:}45$, $60{:}40$, $70{:}30$, and $80{:}20$ and the Newton method.}
    \label{fig:Ex1_NewtonJV}
\end{figure}

\begin{figure}[htp!]
    \centering
    \includegraphics[width=1\linewidth]{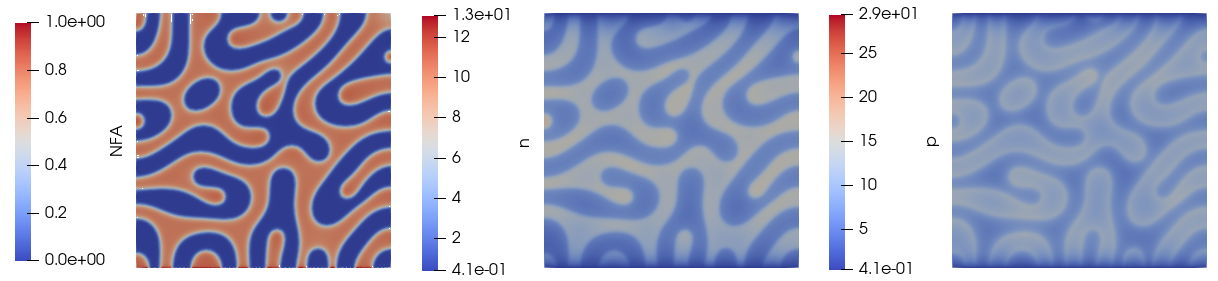}
    
    \hfill
    \includegraphics[width=0.67\linewidth]{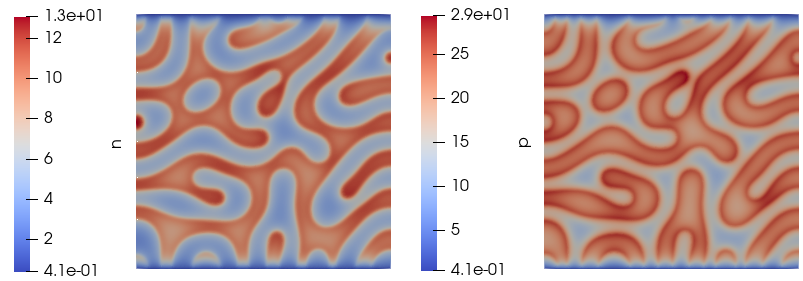}
    \caption{Morphology (left) and snapshots of the electron density (middle) and hole density (right) under forward bias voltages of $0.0,\mathrm{V}$ (top) and $10.0,\mathrm{V}$ (bottom) in 2D, computed using the Gummel method.}
    \label{fig:Ex1_gummnp}
\end{figure}
\begin{figure}[htp!]
    \centering
    \includegraphics[width=0.9\linewidth]{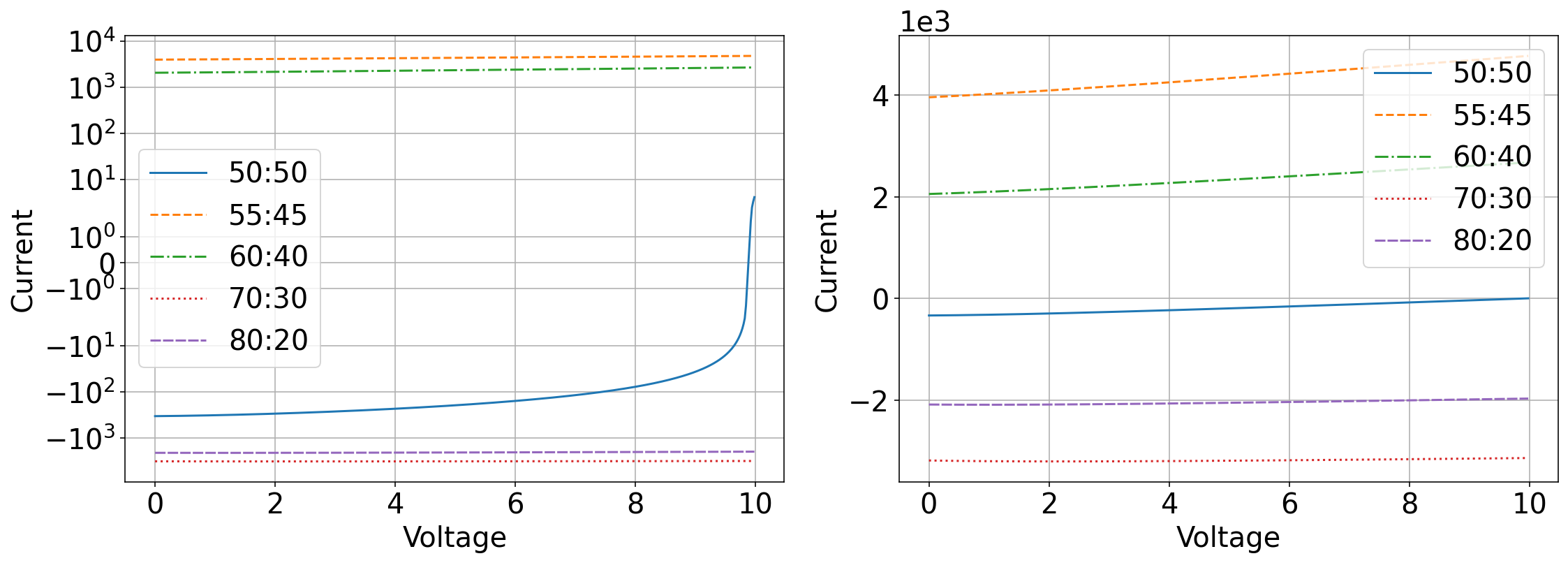}
    \caption{Total currents flowing out of the boundary $\Gamma_{\mathrm{bot/org}}$: semilogarithmic plot (left) and Cartesian plot (right), computed using 2D morphologies with different blend ratios $50{:}50$, $55{:}45$, $60{:}40$, $70{:}30$, and $80{:}20$ with the Gummel method.}
    \label{fig:Ex1_GummelJV}
\end{figure}

Figure~\ref{fig:Ex1_Newtonnp} shows the results of 2D simulations at time $t =0.001$ obtained using the Newton method. The morphology corresponds to a blend ratio of 50:50. As expected, electrons concentrate in the NFA--rich regions while holes accumulate in the polymer--rich regions. The electrical potential and exciton density similarly reflect the underlying morphology. Figure~\ref{fig:Ex1_NewtonJV} presents the current--voltage (IV) characteristics obtained using different 2D morphology blends. As expected, the current increases exponentially with increasing voltage for all blend ratios.

\begin{figure}[htp!]
    \centering
    \includegraphics[width=1\linewidth]{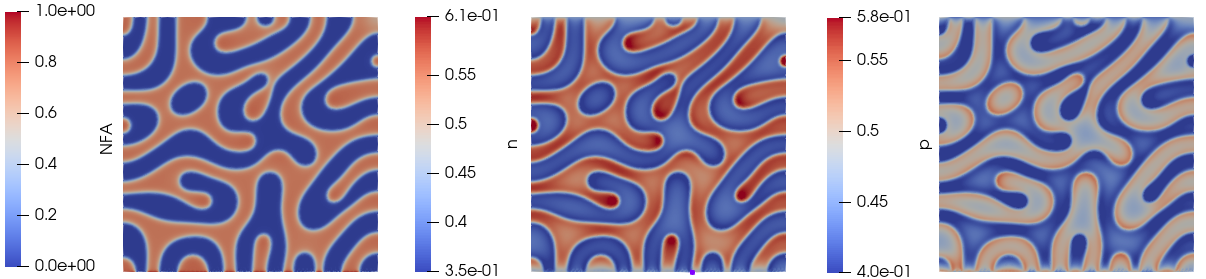}
    \includegraphics[width=1\linewidth]{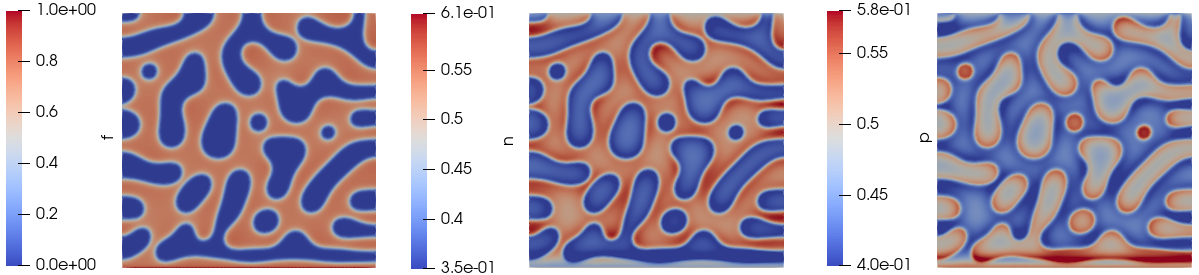}
    \includegraphics[width=1\linewidth]{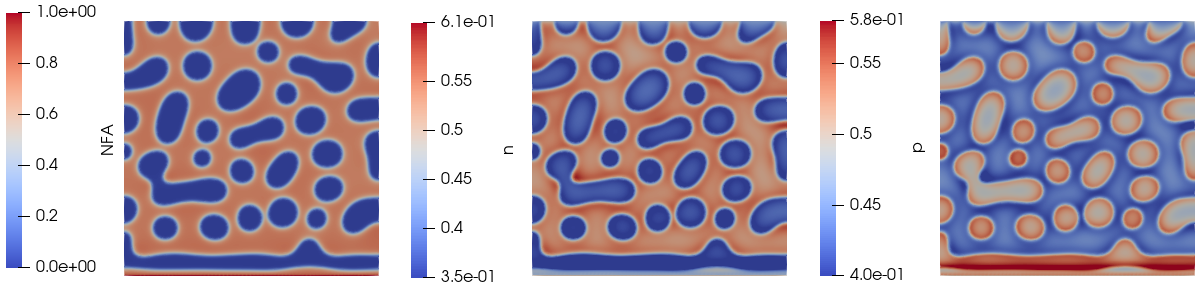}
    \includegraphics[width=1\linewidth]{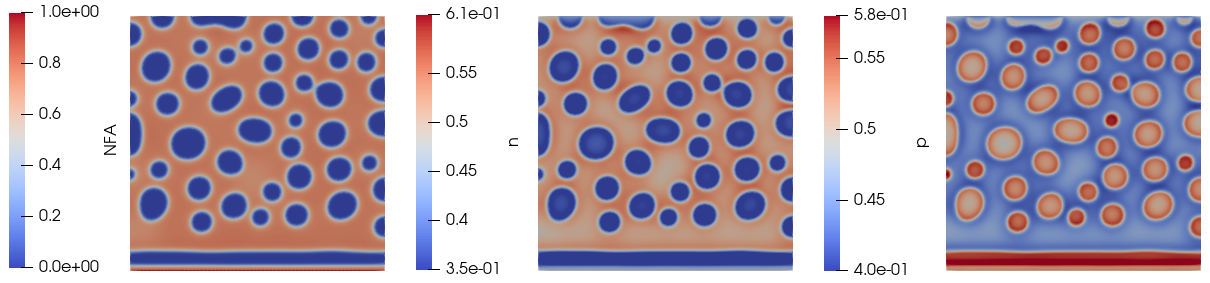}
    \includegraphics[width=1\linewidth]{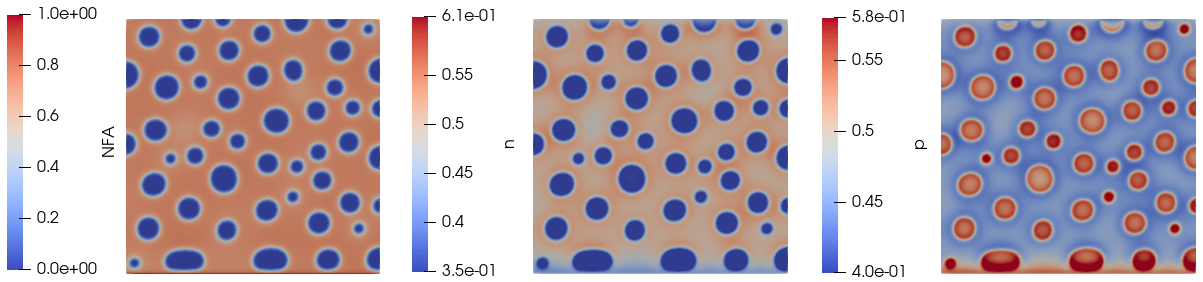}
    \caption{Morphology (left) and snapshots of the electron density (middle) and hole density (right) at time $t = 0.001$ in 2D, computed using the Semi--Newton--Gummel method. From top to bottom, the results correspond to blend ratios $50{:}50$, $55{:}45$, $60{:}40$, $70{:}30$, and $80{:}20$.}
    \label{fig:Ex1_seminp}
\end{figure}



\begin{figure}[htp!]
    \centering
    \includegraphics[width=0.9\linewidth]{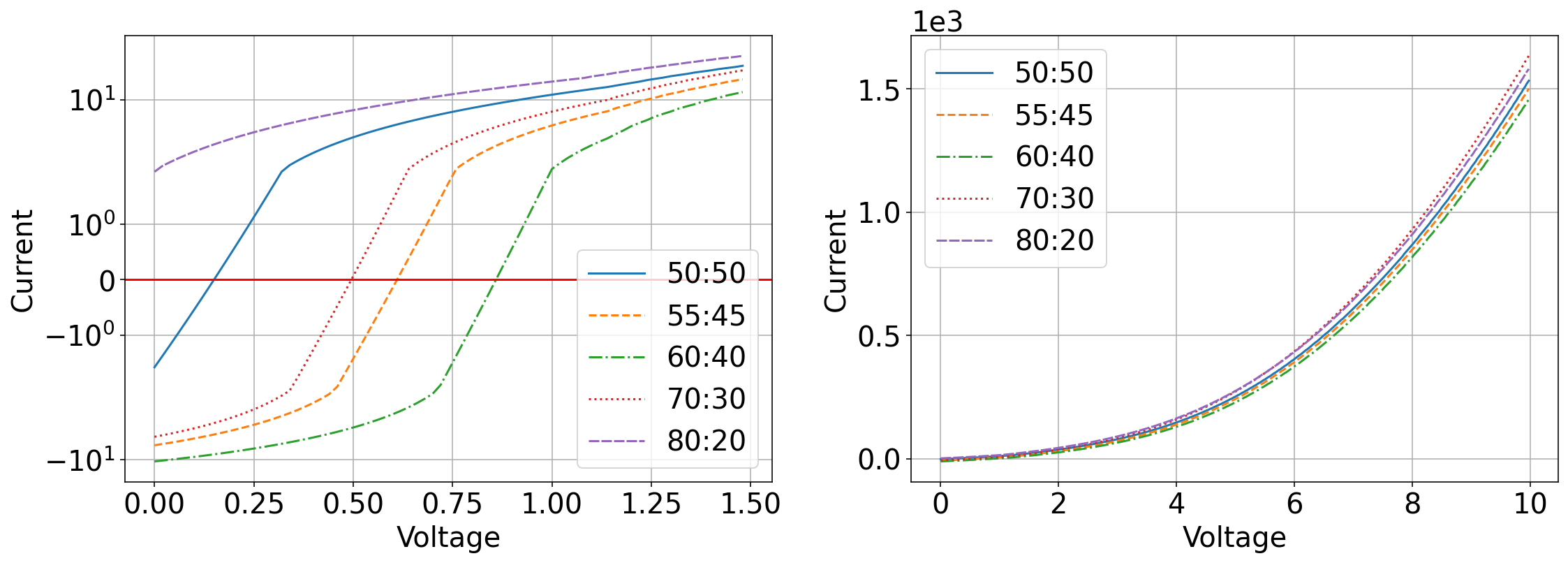}
    \caption{Total currents flowing out of the boundary $\Gamma_{\mathrm{bot/org}}$: semilogarithmic plot (left) and Cartesian plot (right), computed using 2D morphologies with different blend ratios $50{:}50$, $55{:}45$, $60{:}40$, $70{:}30$, and $80{:}20$ and the Semi--Newton--Gummel method.}
    \label{fig:Ex1_JV}
\end{figure}

\begin{table}[htp!]
\centering
\caption{Terminal currents flowing out of $\Gamma_{\mathrm{top/org}}$ and $\Gamma_{\mathrm{bot/org}}$ for different applied external voltages $V_{\mathrm{top}}$ in 2D, computed using the Semi--Newton--Gummel method.}
\begin{footnotesize}
\begin{tabular}{c c cccccc}
\hline
   & $V_{\mathrm{top}}$ & 0.0 & 1.0 & 2.0 & 4.0 & 10.0 \\
\hline
\multirow{2}{*}{50:50} & $I_{\Gamma_{\mathrm{top/org}}}$  & 1.5847e+00 & -1.1236e+01 & -3.7464e+01 & -1.4775e+02 & -1.5407e+03 \\
  & $I_{\Gamma_{\mathrm{bot/org}}}$  & -1.5846e+00 & 1.1236e+01 & 3.7464e+01 & 1.4775e+02 & 1.5407e+03 \\
\hline
\multirow{2}{*}{55:45} & $I_{\Gamma_{\mathrm{top/org}}}$  & 7.1623e+00 & -5.5109e+00 & -3.1237e+01 & -1.3928e+02 & -1.5087e+03 \\
  & $I_{\Gamma_{\mathrm{bot/org}}}$  & -7.1623e+00 & 5.5109e+00 & 3.1237e+01 & 1.3928e+02 & 1.5087e+03 \\
\hline
\multirow{2}{*}{60:40} & $I_{\Gamma_{\mathrm{top/org}}}$  & 1.0413e+01 & -2.0099e+00 & -2.6370e+01 & -1.2995e+02 & -1.4654e+03 \\
  & $I_{\Gamma_{\mathrm{bot/org}}}$  & -1.0413e+01 & 2.0100e+00 & 2.6370e+01 & 1.2996e+02 & 1.4654e+03 \\
\hline
\multirow{2}{*}{70:30} & $I_{\Gamma_{\mathrm{top/org}}}$  & 5.8887e+00 & -7.5930e+00 & -3.7365e+01 & -1.5903e+02 & -1.6440e+03 \\
  & $I_{\Gamma_{\mathrm{bot/org}}}$  & -5.8887e+00 & 7.5931e+00 & 3.7365e+01 & 1.5904e+02 & 1.6440e+03 \\
\hline
\multirow{2}{*}{80:20} & $I_{\Gamma_{\mathrm{top/org}}}$  & -1.9545e+00 & -1.5228e+01 & -4.4690e+01 & -1.6403e+02 & -1.5902e+03 \\
  & $I_{\Gamma_{\mathrm{bot/org}}}$  & 1.9544e+00 & 1.5228e+01 & 4.4690e+01 & 1.6403e+02 & 1.5903e+03 \\
\hline
\end{tabular}
\label{tab:Ex1Current}
\end{footnotesize}
\end{table}

In Figure~\ref{fig:Ex1_gummnp}, the morphology with a blend ratio of 50:50 is analyzed using the Gummel method. Electron and hole densities are plotted for different applied voltages. With increasing external voltage, the electron and hole densities become more localized in the NFA and polymer domains, respectively. However, the I-V curves obtained for different blend ratios do not exhibit an ideal behavior using this method in Figure~\ref{fig:Ex1_GummelJV}. Therefore, as discussed in the previous sections, the results obtained using the hybrid Semi--Newton--Gummel method, which combines the advantages of both approaches, are presented in Figures~\ref{fig:Ex1_seminp} and \ref{fig:Ex1_JV} for different morphologies. For all blend ratios, the electron and hole densities exhibit the expected behavior. Here, it is also clearly observed that the mathematical model accurately captures the actual device behavior. In Figure~\ref{fig:Ex1_JV}, the I-V curves show an exponential increase for all blend ratios. Table~\ref{tab:Ex1Current} reports the terminal currents at the top and bottom boundaries of the organic domain under applied external voltages to verify that our methods guarantee the conservation of the computed terminal currents. The negative sign indicates current inflow. Table~\ref{tab:Ex1Current} shows that terminal currents are conserved; however, floating-point rounding mistakes may cause slight variations in the currents at the two boundaries. Another observation is that at low voltages, the 80:20 blend ratio produces higher current, whereas at higher voltages, the 70:30 blend ratio results in higher current.

\begin{figure}[htp!]
    \centering
    \includegraphics[width=1\linewidth]{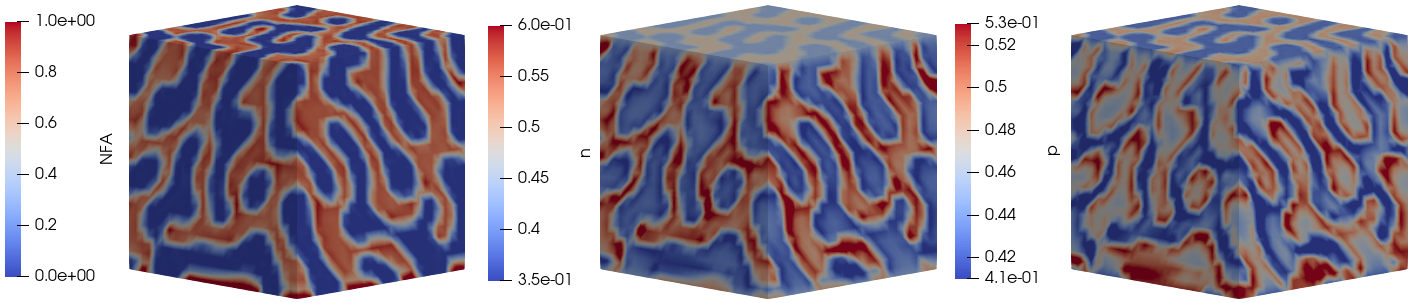}
    \includegraphics[width=1\linewidth]{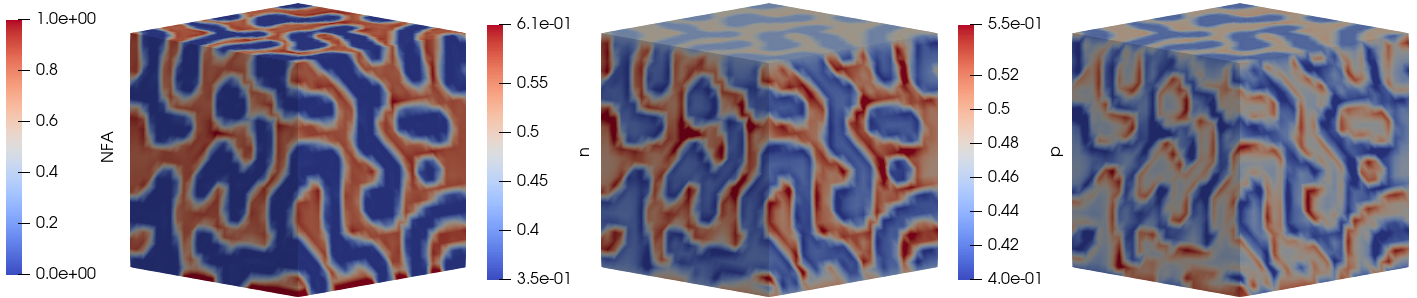}
    \includegraphics[width=1\linewidth]{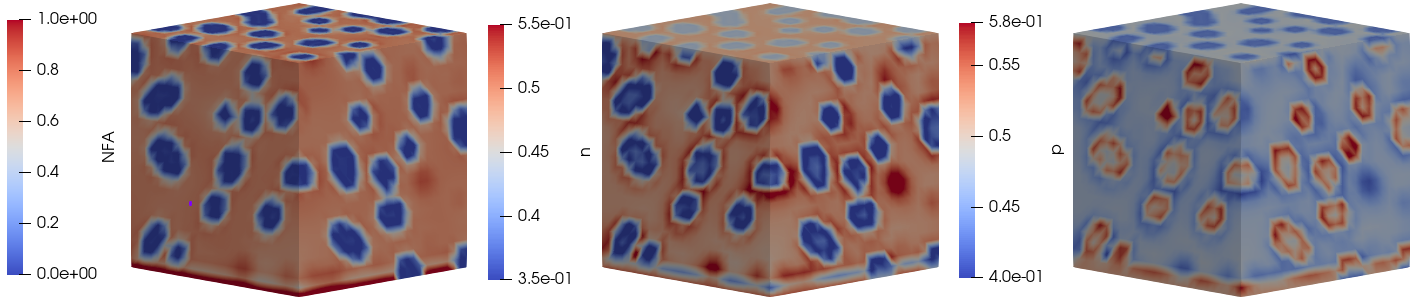}
    \includegraphics[width=1\linewidth]{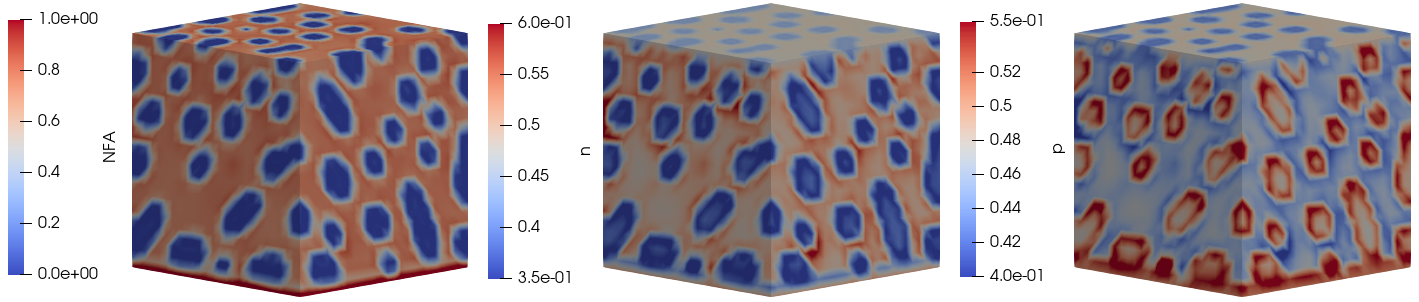}
    \caption{Morphology (left) and snapshots of the electron density (middle) and hole density (right) at time $t = 0.05$ in 3D, computed using the Semi--Newton--Gummel method. From top to bottom, the results correspond to blend ratios $50{:}50$, $55{:}45$, $70{:}30$, and $80{:}20$.}
    \label{fig:Ex1_3D_np}
\end{figure}

\begin{figure}[htp!]
    \centering
    \includegraphics[width=1\linewidth]{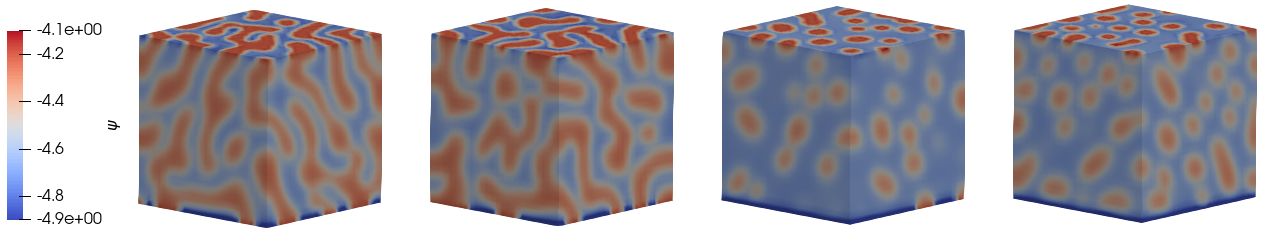}
    \includegraphics[width=1\linewidth]{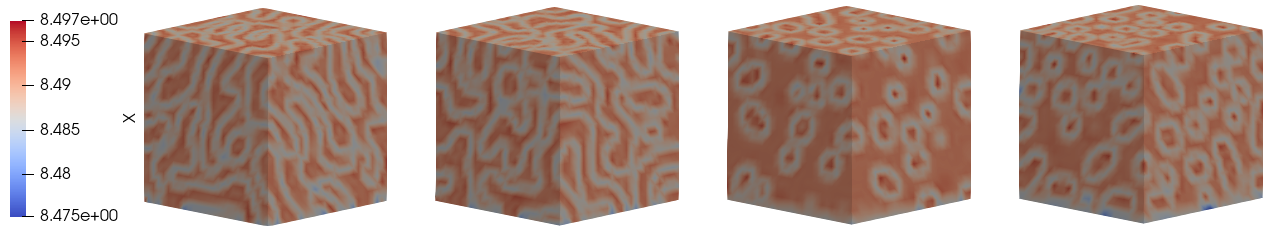}
    \caption{Snapshots of the electrical potential (top) and exciton density (bottom) at time $t = 0.05$ in 3D, computed using the Semi--Newton--Gummel method. From left to right, the results correspond to blend ratios $50{:}50$, $55{:}45$, $70{:}30$, and $80{:}20$.}
    \label{fig:Ex1_3D_psiX}
\end{figure}

\begin{figure}[htp!]
    \centering
    \includegraphics[width=0.9\linewidth]{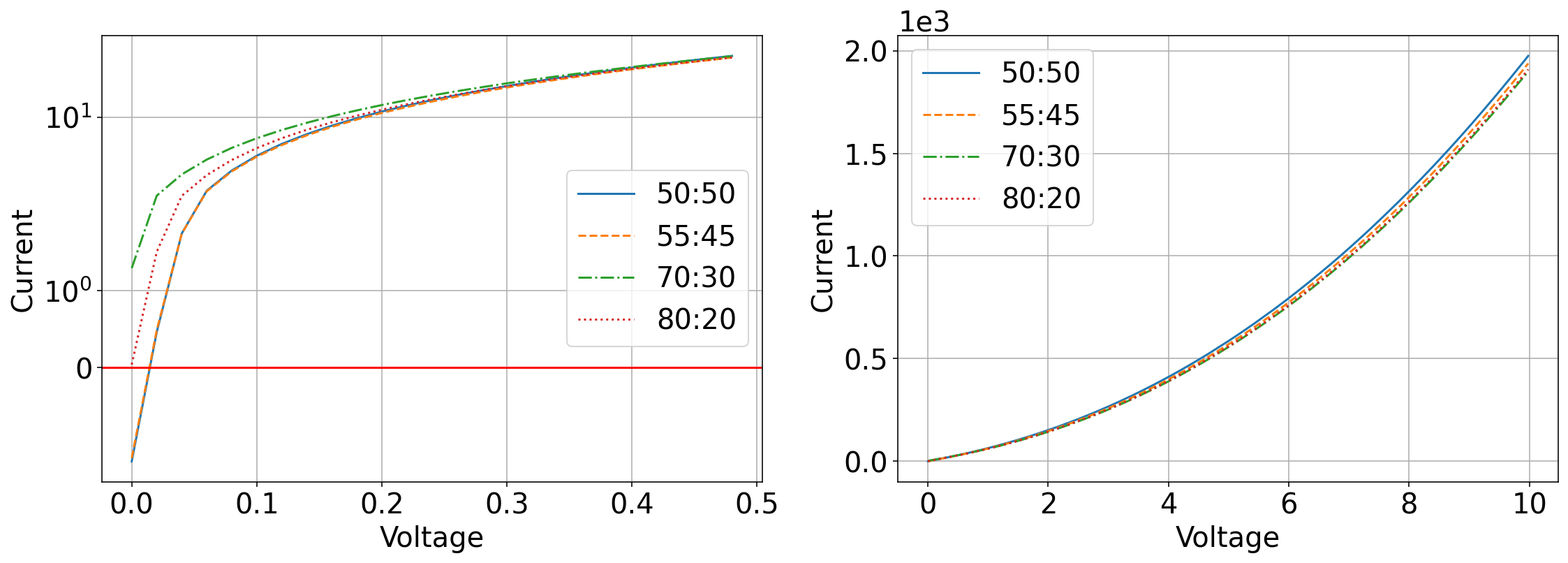}
    \caption{Total currents flowing out of the boundary $\Gamma_{\mathrm{bot/org}}$: semilogarithmic plot (left) and Cartesian plot (right), computed using 3D morphologies with blend ratios $50{:}50$, $55{:}45$, $70{:}30$, and $80{:}20$ and the Semi--Newton--Gummel method. }
    \label{fig:Ex1_3D_JV}
\end{figure}

To further extend our model, we also analyzed 3D morphologies for different blend ratios using the Semi--Newton-–Gummel method, and the results are shown in Figures~\ref{fig:Ex1_3D_np},~\ref{fig:Ex1_3D_psiX}, and ~\ref{fig:Ex1_3D_JV}. Similar to the 2D cases, the electron and hole densities, electrical potential, and exciton density are shaped according to the morphology of the organic domain. From Figure~\ref{fig:Ex1_3D_JV} and Table~\ref{tab:Ex1_current_3D}, it is observed that the current increases exponentially with applied voltage, and at high voltages, the 50:50 blend ratio generates the highest current in the 3D morphologies. Furthermore, conservation of the computed terminal currents is also guaranteed in 3D.  

Finally, to demonstrate the efficiency of the methods, we compared the performance of solvers in this example. Table~\ref{tab:Ex1iterations} reports the maximum number of Newton iterations and \texttt{GMRES} iterations for both 2D and 3D simulations throughout the process. Although the Newton method requires relatively few iterations, it suffers from computational time issues. The Semi--Newton--Gummel method, by leveraging the Gummel decoupling approach instead of solving a large coupled linear system, maintains the maximum iteration numbers within a small range and achieves faster and robust results.

\begin{table}[htp!]
\centering
\caption{Terminal currents flowing out of $\Gamma_{\mathrm{top/org}}$ and $\Gamma_{\mathrm{bot/org}}$ for different external voltages $V_{\mathrm{top}}$ in 3D, computed using the Semi--Newton--Gummel method.}
\begin{footnotesize}
\begin{tabular}{c c c c c c c }
\hline
 & $V_{\mathrm{top}}$ & 0.0 & 1.0 & 2.0 & 4.0 & 10.0 \\
\hline
\multirow{2}{*}{50:50} & $I_{\Gamma_{\mathrm{top/org}}}$ & 1.2243e+00 & -6.3659e+01 & -1.5152e+02 & -4.1063e+02 & -1.9826e+03 \\
 & $I_{\Gamma_{\mathrm{bot/org}}}$ & -1.2299e+00 & 6.3646e+01 & 1.5150e+02 & 4.1059e+02 & 1.9824e+03 \\
\hline
\multirow{2}{*}{55:45} & $I_{\Gamma_{\mathrm{top/org}}}$ & 1.1822e+00 & -6.1763e+01 & -1.4680e+02 & -3.9873e+02 & -1.9452e+03 \\
 & $I_{\Gamma_{\mathrm{bot/org}}}$ & -1.1873e+00 & 6.1751e+01 & 1.4678e+02 & 3.9869e+02 & 1.9450e+03 \\
\hline
\multirow{2}{*}{70:30} & $I_{\Gamma_{\mathrm{top/org}}}$ & -1.2957e+00 & -6.0911e+01 & -1.4312e+02 & -3.8893e+02 & -1.9068e+03 \\
 & $I_{\Gamma_{\mathrm{bot/org}}}$ & 1.2898e+00 & 6.0898e+01 & 1.4310e+02 & 3.8888e+02 & 1.9067e+03 \\
\hline
\multirow{2}{*}{80:20} & $I_{\Gamma_{\mathrm{top/org}}}$ & -3.6972e-02 & -6.0854e+01 & -1.4390e+02 & -3.9103e+02 & -1.9158e+03 \\
 & $I_{\Gamma_{\mathrm{bot/org}}}$ & 3.2300e-02 & 6.0842e+01 & 1.4388e+02 & 3.9099e+02 & 1.9157e+03 \\
\hline
\end{tabular}
\label{tab:Ex1_current_3D}
\end{footnotesize}
\end{table}

\begin{table}[h!]
\centering
\caption{Maximum Newton and \texttt{GMRES} iterations for different methods in 2D (and 3D).}
\begin{tabular}{llcc}
\hline
 & Equation & Newton Iteration & \texttt{GMRES} Iteration \\
\hline
\multirow{2}{*}{Newton} & Poisson-Elec-Hole & 8 (8) & 5 (4) \\
 & Exciton & $ - $ & 5 (4) \\
\hline
\multirow{4}{*}{Gummel} & Poisson & $ - $  & 12 (12) \\
 & Electron & $ - $  & 11 (10) \\
 & Hole & $ - $  & 10 (10) \\
 & Exciton & $ - $  & 4 (4) \\
\hline
\multirow{4}{*}{\makecell{Semi--Newton \\ Gummel}} & Poisson & 8 (8) & 4 (3) \\
 & Electron & 5 (4) & 4 (6) \\
 & Hole & 5 (4) & 4 (6)\\
 & Exciton & $ - $  & 5 (5)\\
\hline
\end{tabular}
\label{tab:Ex1iterations}
\end{table}

\subsection{Example 2:} \label{sec_example2}
In this subsection, we present the result using the morphologies obtained by the four-component evaporative phase-field model. In generating the morphologies, we initialized the domain with a size of $10\times10\times10$ in 3D ($10\times10$ in 2D), and the blend ratios 50:50, where the first $h=8$ is occupied by a liquid mixture of polymer, non--fullerene acceptor (NFA), and solvent, and the remaining $h=2$ represents the gas phase. 
\begin{figure}[htp!]
    \centering
    \includegraphics[width=1\linewidth]{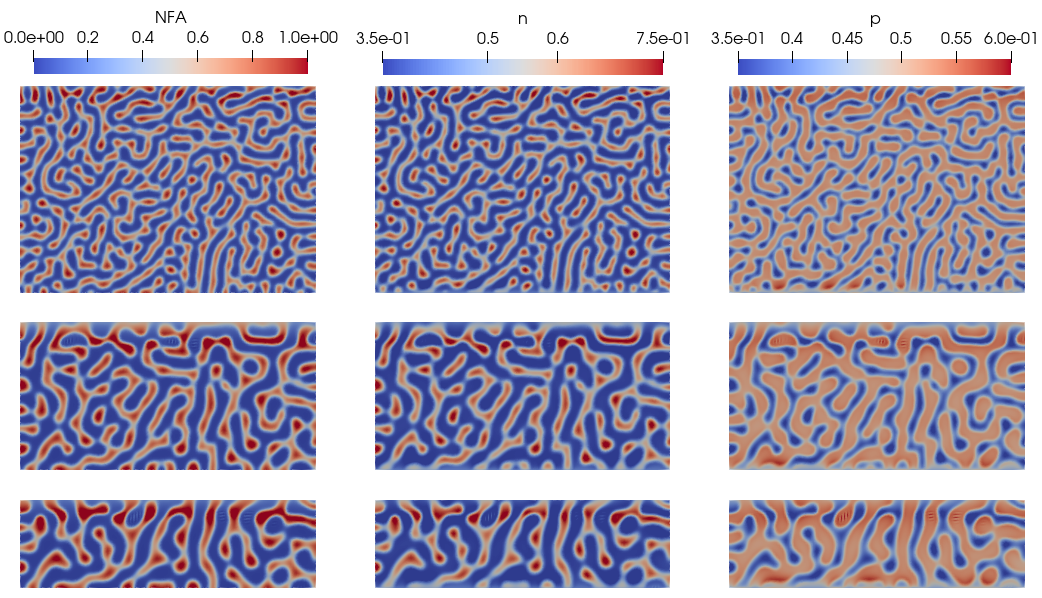}
    \caption{Morphology (left) and snapshots of the electron density (middle) and hole density (right) at time $t = 0.05$ in 2D, computed using the Semi--Newton--Gummel method. From top to bottom, the results correspond to film heights $h = 7, 5,$ and $3$.}
    \label{fig:Ex2_2D_np}
\end{figure}

\begin{figure}[htp!]
    \centering
    \includegraphics[width=1\linewidth]{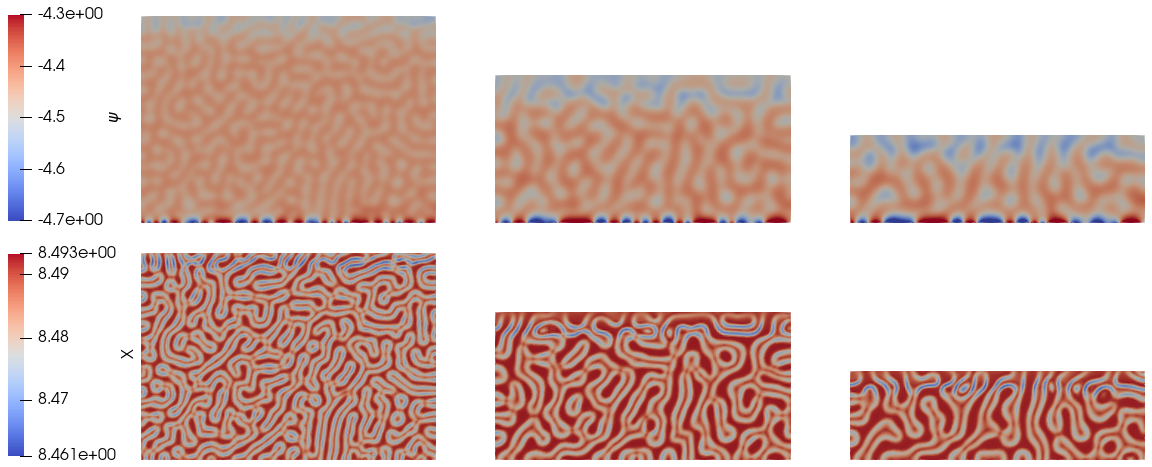}
    \caption{Snapshots of the electrical potential (top) and the exciton density (bottom) at time $t = 0.05$ in 2D, computed using the Semi--Newton--Gummel method. From top to bottom, the results correspond to film heights $h = 7, 5,$ and $3$.}
    \label{fig:Ex2_2D_psiX}
\end{figure}

We use $n_x = 200$ and $n_y = 200$ grid points in 2D, and $n_x = 50$, $n_y = 50$, and $n_z = 50$ grid points in 3D. Morphologies corresponding to different film heights will be analyzed to test the performance of the model.

\begin{figure}[htp!]
    \centering
    \includegraphics[width=0.9\linewidth]{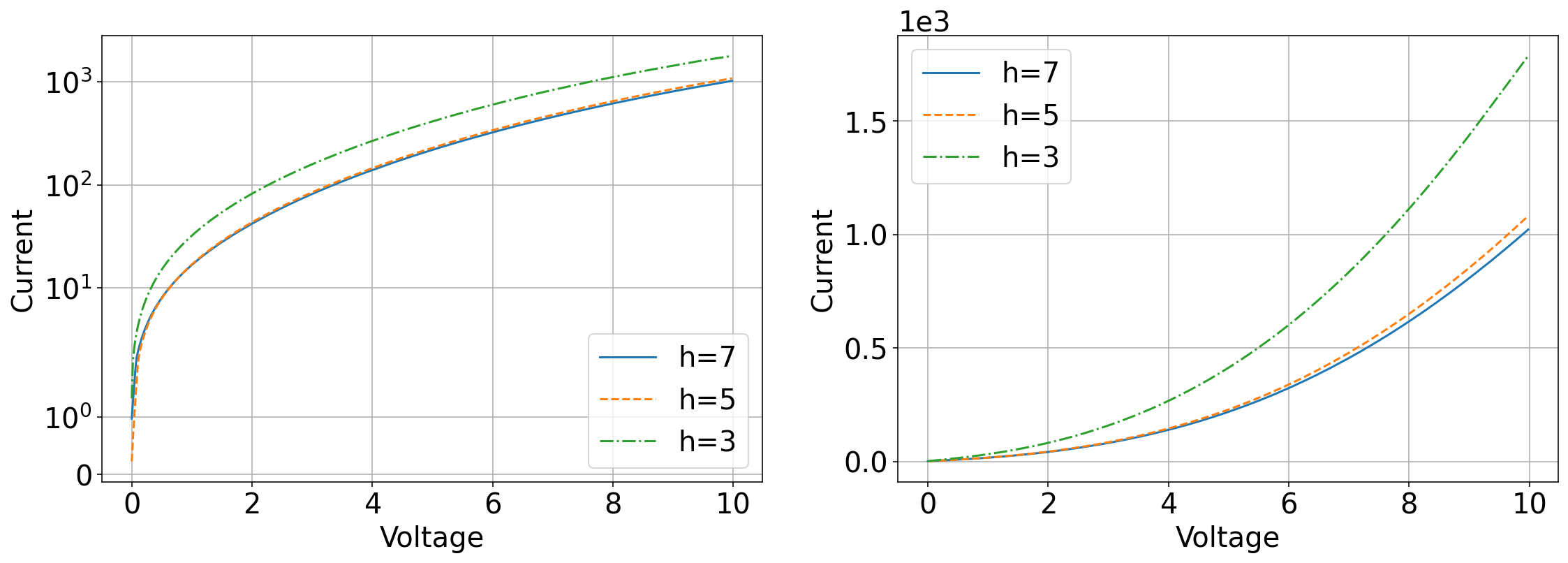}
    \caption{Total currents flowing out of the boundary $\Gamma_{\mathrm{bot/org}}$: semilogarithmic plot (left) and Cartesian plot (right), computed using 2D morphologies with different film heights $h = 7, 5,$ and $3$ and the Semi--Newton--Gummel method.}
    \label{fig:Ex2_2D_JV}
\end{figure}

In these morphologies, we observe that as the solvent evaporates, the film thickness decreases and more separated donor–acceptor regions emerge, leading to a coarser structure, as illustrated in Figure \ref{fig:Ex2_2D_np}. Using 2D morphologies with different film heights, the results obtained with the Semi--Newton--Gummel method are presented in Figures \ref{fig:Ex2_2D_np} and \ref{fig:Ex2_2D_psiX}. Similar to the previous example, despite employing different film thicknesses, the expected behavior of the electron and hole densities, electrical potential, and exciton density is successfully captured. From the I-V curves shown in Figure \ref{fig:Ex2_2D_JV} and the results reported in Table \ref{tab:Ex2_2D}, it is observed that the current increases exponentially with increasing applied voltage. In addition, thinner films produce higher currents. Despite the variation in organic domain heights, conservation of the computed terminal currents at the top and bottom boundaries is maintained.

\begin{table}[htp!]
\centering
\caption{Terminal currents flowing out of $\Gamma_{\mathrm{top/org}}$ and $\Gamma_{\mathrm{bot/org}}$ for different external voltages $V_{\mathrm{top}}$ in 2D.}
\begin{footnotesize}
\begin{tabular}{c c c c c c c }
\hline
 & $V_{\mathrm{top}}$ & 0.0 & 1.0 & 2.0 & 4.0 & 10.0 \\
\hline
\multirow{2}{*}{h=7} & $I_{\Gamma_{\mathrm{top/org}}}$ & -9.7356e-01 & -1.6708e+01 & -4.1998e+01 & -1.3937e+02 & -1.0267e+03 \\
 & $I_{\Gamma_{\mathrm{bot/org}}}$ & 9.7358e-01 & 1.6708e+01 & 4.1998e+01 & 1.3937e+02 & 1.0267e+03 \\
\hline
\multirow{2}{*}{h=5} & $I_{\Gamma_{\mathrm{top/org}}}$ & -2.2554e-01 & -1.6894e+01 & -4.3353e+01 & -1.4543e+02 & -1.0872e+03 \\
 & $I_{\Gamma_{\mathrm{bot/org}}}$ & 2.2554e-01 & 1.6894e+01 & 4.3353e+01 & 1.4543e+02 & 1.0872e+03 \\
\hline
\multirow{2}{*}{h=3} & $I_{\Gamma_{\mathrm{top/org}}}$ & -1.3310e+00 & -3.2181e+01 & -8.2152e+01 & -2.6722e+02 & -1.7944e+03 \\
 & $I_{\Gamma_{\mathrm{bot/org}}}$ & 1.3311e+00 & 3.2181e+01 & 8.2152e+01 & 2.6722e+02 & 1.7944e+03 \\
\hline
\end{tabular}
\label{tab:Ex2_2D}
\end{footnotesize}
\end{table}

\begin{figure}[htp!]
    \centering
    \includegraphics[width=1\linewidth]{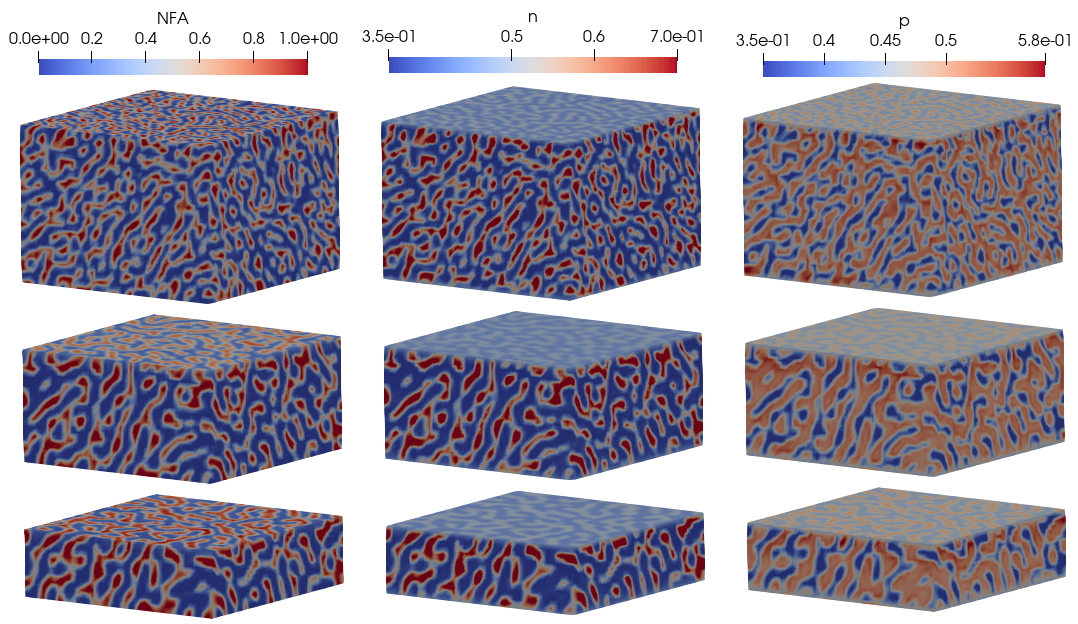}
    \caption{Morphology (left) and snapshots of the electron density (middle) and hole density (right) at time $t = 0.05$ in 3D, computed using the Semi--Newton--Gummel method. From top to bottom, the results correspond to film heights $h = 7, 5,$ and $3$.}
    \label{fig:Ex2_3D_np}
\end{figure}

\begin{figure}[htp!]
    \centering
    \includegraphics[width=1\linewidth]{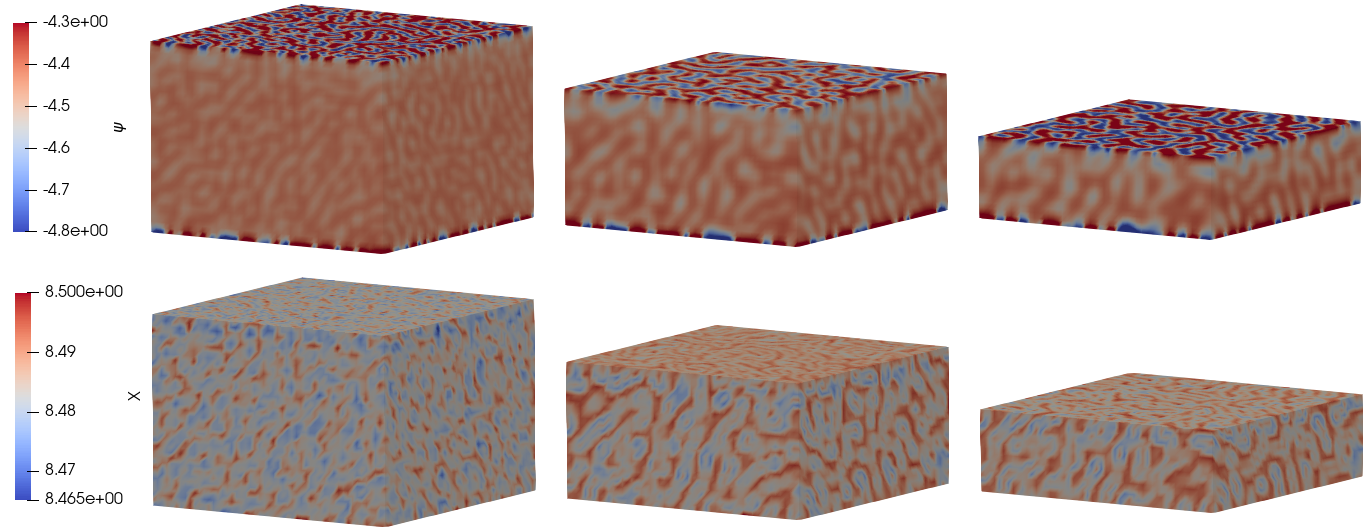}
    \caption{Snapshots of the electrical potential (top) and the exciton density (bottom) at time $t = 0.05$ in 3D, computed using the Semi--Newton--Gummel method. From top to bottom, the results correspond to film heights $h = 7, 5,$ and $3$.}
    \label{fig:Ex2_3D_psiX}
\end{figure}

\begin{figure}[htp!]
    \centering
    \includegraphics[width=0.9\linewidth]{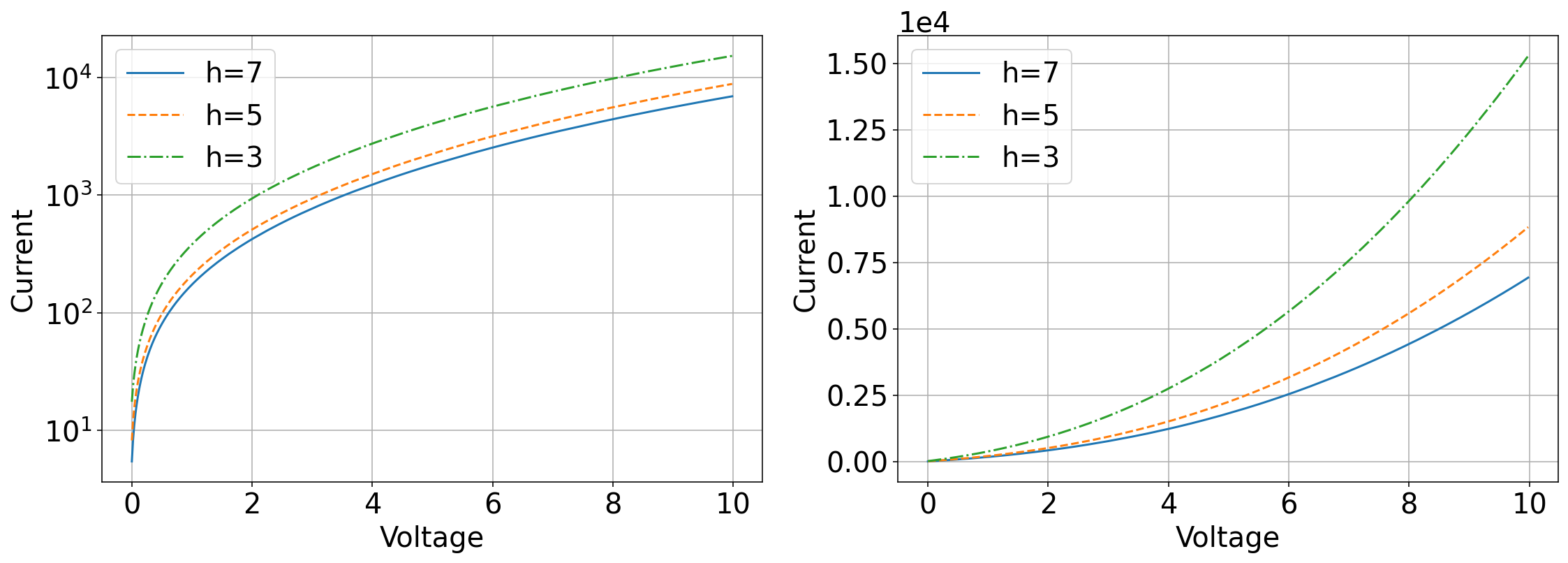}
    \caption{Total currents flowing out of the boundary $\Gamma_{\mathrm{bot/org}}$: semilogarithmic plot (left) and Cartesian plot (right), computed using 3D morphologies with different film heights $h = 7, 5,$ and $3$ and the Semi--Newton--Gummel method.}
    \label{fig:Ex2_3D_JV}
\end{figure}

Following the same approach, this example is also extended to the 3D model. As shown in Figures \ref{fig:Ex2_3D_np} and \ref{fig:Ex2_3D_psiX}, the mathematical model reproduces the expected physical behavior: electrons concentrate in the NFA--rich regions, while holes accumulate in the polymer--rich regions. Moreover,  thinner films generate higher currents, and when comparing the 2D I-V curves, a clear separation between the currents corresponding to different film heights is observed, see Figure \ref{fig:Ex2_3D_JV} and Table \ref{tab:Ex2_3D}.

Finally, the performance of the 2D and 3D solvers is compared for different film heights in Table \ref{tab:Ex2_iter}. As in Example \ref{sec_example1}, the Semi--Newton--Gummel method keeps the maximum numbers of Newton and GMRES iterations within a small range, resulting in robust and efficient performance for the different film heights.

\begin{table}[htp!]
\centering
\caption{Terminal currents flowing out of $\Gamma_{\mathrm{top/org}}$ and $\Gamma_{\mathrm{bot/org}}$ for different external voltages $V_{\mathrm{top}}$ in 3D.}
\begin{footnotesize}
\begin{tabular}{c c c c c c c }
\hline
 & $V_{\mathrm{top}}$ & 0.0 & 1.0 & 2.0 & 4.0 & 10.0 \\
\hline
\multirow{2}{*}{h=7} & $I_{\Gamma_{\mathrm{top/org}}}$ & -5.4540e+00 & -1.7396e+02 & -4.2361e+02 & -1.2316e+03 & -6.9571e+03 \\
 & $I_{\Gamma_{\mathrm{bot/org}}}$ & 5.4542e+00 & 1.7396e+02 & 4.2361e+02 & 1.2316e+03 & 6.9571e+03 \\
\hline
\multirow{2}{*}{h=5} & $I_{\Gamma_{\mathrm{top/org}}}$ & -8.2400e+00 & -2.0909e+02 & -5.1089e+02 & -1.5125e+03 & -8.8669e+03 \\
 & $I_{\Gamma_{\mathrm{bot/org}}}$ & 8.2328e+00 & 2.0907e+02 & 5.1086e+02 & 1.5124e+03 & 8.8665e+03 \\
\hline
\multirow{2}{*}{h=3} & $I_{\Gamma_{\mathrm{top/org}}}$ & -1.7532e+01 & -3.8186e+02 & -9.3873e+02 & -2.7498e+03 & -1.5346e+04 \\
 & $I_{\Gamma_{\mathrm{bot/org}}}$ & 1.7529e+01 & 3.8185e+02 & 9.3872e+02 & 2.7498e+03 & 1.5345e+04 \\
\hline
\end{tabular}
\label{tab:Ex2_3D}
\end{footnotesize}
\end{table}

\begin{table}[h!]
\centering
\caption{Maximum Newton and \texttt{GMRES} iterations for different film heights $h = 7, 5,$ and $3$ in 2D (and 3D).}
\begin{footnotesize}
\begin{tabular}{llccc ccc}
\hline
 & Equation 
 & \multicolumn{3}{c}{Newton Iteration}
 & \multicolumn{3}{c}{\texttt{GMRES} Iteration} \\
\cline{3-5} \cline{6-8}
 & 
 & $h=7$ & $h=5$ & $h=3$
 & $h=7$ & $h=5$ & $h=3$ \\
\hline
\multirow{2}{*}{Newton} 
 & Poisson-Elec-Hole 
 & 4 (4) & 4 (4) & 4 (4) 
 & 4 (3) & 4 (3) & 4 (3)  \\
 & Exciton 
 & $-$ & 5 (6) & 5 (6) 
 & 5 (6) & 5 (6) & 5 (6) \\
\hline
\multirow{4}{*}{Gummel} 
 & Poisson 
 & $-$ & 11 (13) & 11 (12) 
 & 11 (12) & 11 (12) & 11 (12) \\
 & Electron 
 & $-$ & 15 (25) & 10 (12) 
 & 10 (12) & 10 (12) & 10 (12) \\
 & Hole 
 & $-$ & 10 (12) & 9 (10) 
 & 9 (10) & 9 (10) & 9 (10) \\
 & Exciton 
 & $-$ & 4 (4) & 4 (4) 
 & 4 (4) & 4 (4) & 4 (4) \\
\hline
\multirow{4}{*}{\makecell{Semi--Newton \\ Gummel}} 
 & Poisson 
 & 8 (8) & 8 (8) & 8 (8) 
 & 4 (3) & 4 (3) & 4 (3) \\
 & Electron 
 & 4 (4) & 4 (4) & 4 (4) 
 & 4 (5) & 4 (5) & 4 (5) \\
 & Hole 
 & 4 (4) & 4 (4) & 4 (4) 
 & 4 (5) & 4 (5) & 4 (5) \\
 & Exciton 
 & $-$ & 5 (5) & 5 (5) 
 & 5 (5) & 5 (5) & 5 (5) \\
\hline
\end{tabular}
\end{footnotesize}
\label{tab:Ex2_iter}
\end{table}

\section{Conclusions}
\label{sec:conclusions}
In this study, we presented a comprehensive in silico framework for investigating the performance of organic solar cells. The mathematical model establishes a direct link between the nanomorphology of the active layer, exciton dynamics, and performance indicators at the device level. Although the derived model makes use of some simplifications, like the omission of any thermal effects or the neglect of the energetic disorder associated with charge trapping, it still offers detailed insight into the influence of morphological variations on charge generation and transport. 

The resulting system of nonlinear partial differential equations poses significant numerical challenges due to its strong coupling, pronounced nonlinearity, and large number of degrees of freedom, especially in three dimensions. To overcome these challenges, we employed a Semi--Newton--Gummel strategy combined with iterative linear solvers and appropriate preconditioning techniques. The numerical results demonstrate that the proposed approach provides a robust, efficient framework for two- and three-dimensional simulations. Specifically, the decoupled Gummel iteration improves numerical stability, and the \texttt{GMRES} method efficiently handles large-scale problems resulting from fine spatial discretizations.

Considering device morphology to play a critical role in OPV performance, the developed pipeline shows to be a promising tool for future production and design choices. To derive additional meaningful conclusions for physically relevant settings, as a next step, parameters must be matched and resulting I-V curves must be compared to experimental measurements.





\section*{Acknowledgments}
We thank the Deutsche Forschungsgemeinschaft (DFG) for funding this work (Research Unit FOR 5387 POPULAR, project no.\ 461909888).

\bibliographystyle{siamplain}
\bibliography{refs}

\end{document}